\documentclass[12pt]{article} 
\pdfoutput=1
\usepackage[a4paper]{geometry}
\usepackage[utf8]{inputenc} \usepackage[english]{babel} \usepackage[T2A]{fontenc}
\usepackage{amsmath,amssymb} \usepackage{euler,eucal}
\usepackage{myconcrete} 
\usepackage[unicode=true,plainpages=false]{hyperref} 
\hypersetup{colorlinks=true,linkcolor=magenta,anchorcolor=magenta,urlcolor=blue,citecolor=blue}
\usepackage{algorithmic} 
\usepackage[ruled]{algorithm} 
\usepackage{graphicx}
\usepackage{setspace}
\def\title#1{\phantom{.}\vskip 5mm
      \begin{center}\begin{doublespace}{ \LARGE\sc  #1}\end{doublespace}\end{center}}

\def\author#1{\begin{center}\large{#1}\end{center}}

\def\address#1{\begin{center}\footnotesize\textit{#1}\end{center}}

\def\date#1{\vskip 3mm\begin{center}\large{#1}\end{center}\vskip 10mm}

\def\affilnum#1{${}^{\@fnsymbol{#1}}$}
\def\affil#1{${}^{\@fnsymbol{#1}}$}

\def\thanks#1{${}^{\@crssymbol{#1}}$}
\def\grant#1{\def\thefootnote{\$}\footnotetext{#1} }

\def\email#1{\texttt{#1}}
\def\keywordsname{Keywords}
\def\aabstractname{Abstract.}
\renewenvironment{abstract}{\small \quotation {\bfseries \aabstractname} }{\endquotation}
\def\keywords#1{\vskip 2mm\par\noindent{\small\emph{\keywordsname:} #1}}
\usepackage{caption}
\usepackage{amsthm}

\theoremstyle{definition}
\newtheorem{remark}{Remark}
\newtheorem{statement}{Statement}
\newtheorem{corollary}{Corollary}
\newtheorem{theorem}{Theorem}
\def\hat{\mathaccent"705E }
\let\eps=\varepsilon 
  
\let\leq=\leqslant \let\geq=\geqslant

\def\ppkeywords{Multidimensional arrays, sparse tensors, structured tensors, Tucker approximation,  Krylov subspace methods, Wedderburn rank reduction, fast compression}

\hypersetup{pdfauthor={S. A. Goreinov, I. V. Oseledets, D. V. Savostyanov},
            pdftitle={Wedderburn rank reduction and Krylov subspace method for tensor approximation. Part 1: Tucker case},
            pdfsubject={Research article on tensor approximations, submitted to SIAM J. Scientific Computing},
            pdfkeywords={\ppkeywords}
            }

\def\cntr#1{\begin{center}\small\it {#1} \end{center}}
\def\trans{T}
\def\t{\trans}

\def\eps{\varepsilon}
\def\phi{\varphi}

\def\A{\mathbf{A}}
\def\B{\mathbf{B}}
\def\C{\mathbf{C}}
\def\F{\mathbf{F}}
\def\G{\mathbf{G}}
\def\H{\mathbf{H}}

\def\O{\mathcal{O}}

\def\vv{\hat v}
\def\ww{\hat w}

\def\yy{\hat y}
\def\zz{\hat z}

\def\x{\mathbin{\times}}
\def\o{\mathbin{\raise1pt\hbox{$\scriptscriptstyle\mathord\otimes$}}}
\def\hadamar{\mathbin{\raise1pt\hbox{$\scriptstyle\mathord\odot$}}}

\def\nrm{\mathop{\mathtt{nrm}}\nolimits}
\def\err{\mathop{\mathtt{err}}\nolimits}

\def\tol{\mathtt{tol}}
\def\rmax{r_{\mathrm{max}}}

\renewcommand{\f}[2]{\left\langle {#1},\:{#2} \right\rangle}

\def\rank{\mathop{\mathrm{rank}}\nolimits}
\def\Span{\mathop{\mathrm{span}}\nolimits}
\def\Kron{\mathop{\mathbf{Kron}}\nolimits}

\def\diag{\mathop{\mathrm{diag}}\nolimits}

\def\Sum{\sum\limits}
\def\eqdef{\mathrel{\stackrel{\mathrm{def}}{=}}}

\def\nnn{{n\times n\times n}}
\def\pals{p_\mathrm{als}}
\def\ppow{p_\mathrm{pow}}
\def\updu{{\mathit{upd}U}}
\def\updv{{\mathit{upd}V}}
\def\updw{{\mathit{upd}W}}
\def\updx{{\mathit{upd}X}}
\def\updy{{\mathit{upd}Y}}
\def\updz{{\mathit{upd}Z}}
\def\xwcp{x^{\{\mathrm{wcp}\}}}
\def\xgk{x^{\{\mathrm{lnc}\}}}
\def\Xwcp{X^{\{\mathrm{wcp}\}}}
\def\Xgk{X^{\{\mathrm{lnc}\}}}
\def\tn{\mathbf{2}}

\begin{document}
\bibliographystyle{hsiam}

\title{Wedderburn rank reduction and \\ Krylov subspace method for tensor approximation. Part 1: Tucker case}
\author{S. A. Goreinov, I. V. Oseledets, D. V. Savostyanov}
\address{Institute of Numerical Mathematics, Russian Academy of Sciences, \\
                     Russia, 119333 Moscow, Gubkina 8 \\ 
                  \email{[sergei.goreinov,ivan.oseledets,dmitry.savostyanov]@gmail.com}
           }
\grant{This work was supported by 
         RFBR grants 08-01-00115, 09-01-00565, 09-01-12058, 10-01-00757, RFBR/DFG grant 09-01-91332,
         Russian Federation Gov. contracts No. $\Pi 1178,$ $\Pi 1112,$ $\Pi 940$ and $14.740.11.0345$, 
         Rus. President grant $\mbox{MK--}127.2009.1$ 
         and  Priority research program of Dep. Math. RAS No. $3$ and $5.$
         Part of this work was done during the stay of the second and third authors in the Max-Planck Institute for Mathematics in Sciences in Leipzig (Germany). 
         }
\date{October 19, 2010}

\begin{abstract}
New algorithms are proposed for the Tucker approximation of a 3-tensor, that access it using only the tensor-by-vector-by-vector multiplication subroutine. 
In the matrix case, Krylov methods are methods of choice to approximate the dominant column and row subspaces of a sparse or structured matrix given through the matrix-by-vector multiplication subroutine. 
Using the Wedderburn rank reduction formula, we propose an algorithm of matrix approximation that computes Krylov subspaces and allows generalization to the tensor case. 
Several variants of proposed tensor algorithms differ by pivoting strategies, overall cost and quality of approximation. 
By convincing numerical experiments we show that the proposed methods are faster and more accurate than the minimal Krylov recursion, proposed recently by Eld\'en and Savas.
\keywords{\ppkeywords}
\par \noindent \emph{AMS classification:} 15A23, 15A69, 65F99
\end{abstract}

\section{Introduction}
In this paper we focus on algorithms for the low-rank approximation of large three-dimensional arrays (tensors),
that play increasingly important role in many applications.
Throughout the paper, a \emph{tensor} $\A=[a_{ijk}]$ means an array with three indices.
The indices are also called \emph{modes} or \emph{axes}. 
The number of allowed values for a mode index is called the \emph{mode size}.
In numerical work with tensors of large mode sizes it is crucial to look for the data sparse structures. 
Most common are the following.

The \emph{canonical decomposition}  \cite{hitchcock-sum-1927,cc-parafac-1970,harshman-parafac-1970} (or \emph{canonical approximation}, if the right-hand side does not give $\A$ exactly) of a tensor $\A=[a_{ijk}]$ reads
$$
\A = \Sum_{s=1}^R u_s \o v_s \o w_s, \qquad   a_{ijk} = \Sum_{s=1}^R u_{is} v_{js} w_{ks}. \eqno{(C)}
$$
The minimal possible number of summands is called the \emph{tensor rank} or the \emph{canonical rank} of a given tensor $\A.$
However, the canonical decomposition/approximation of a tensor with minimal value of $R$ is an ill-posed and computationally unstable problem \cite{desilva-2008}, and among several practical algorithms none is known to be absolutely robust. 

The (truncated
) \emph{Tucker decomposition/approximation}  \cite{Tucker} of $\A$ reads 
$$
\A = \G \x_1 U \x_2 V \x_3 W, \qquad 
      a_{ijk} = \Sum_{p=1}^{r_1}\Sum_{q=1}^{r_2} \Sum_{s=1}^{r_3} g_{pqs} u_{ip} v_{jq} w_{ks}. \eqno{(T)}
$$
The quantities $r_1,r_2,r_3$ are referred to as the \emph{Tucker ranks} or the \emph{mode ranks},
the $r_1 \x r_2 \x r_3$ tensor $\G=[g_{pqs}]$ is called the \emph{Tucker core}.
In $d$ dimensions, the memory to store the $r\x r\x \ldots \x r$ core is $r^d,$ that is usually beyond affordable for large $d$ and modest $r$ (so-called \emph{curse of dimensionality}).
In three dimensions for $r\sim 100$ the storage for the core is small and the Tucker decomposition can be used efficiently.

Here and below, the symbol $\x_l$ designates the multiplication of a tensor by a matrix along the $l$-th mode. 
For example, $\B = \A\x_2 M$ means summation on $2^\mathrm{nd}$ index $b_{ijk} = \sum_{j'} m_{jj'} a_{ij'k}.$ 
The notation for tensor operations is not yet standard (for current state of art see review \cite{kolda-review-2009}),
we use the one proposed by de Lathauwer in the article on \emph{multilinear SVD} \cite{lathauwer-svd-2000} (or higher-order SVD, HOSVD). 
In~\cite{lathauwer-svd-2000} the Tucker format arises with additional constraints: orthogonality of factors, that is also assumed in our paper,  and all-orthogonality of the core, that we will relax. 
For a tensor given as a full array of elements, the multilinear SVD provides a quasi-optimal Tucker approximation which further can be refined by different iterative methods such that Tucker-ALS~\cite{tuckerals-1980,lathauwer-rank1-2000}, Newton-Grassmann~\cite{eldensavas-2009,lathauwer-difnew-2009}, etc. 
It reduces to the SVD for three \emph{unfolding matrices} (see~\cite{lathauwer-svd-2000} and~\eqref{eq:unf} later in this paper) and costs $\O(n^4)$ operations for $n_1 \x n_2 \x n_3$ tensor.%
\footnote{We always assume $n_1=n_2=n_3=n$ and $r_1=r_2=r_3=r$ in complexity estimates}
It is clearly too much for large $n,$ and we should look for alternatives. 

In the matrix case two classes of fast methods are popular: cross algorithms and Krylov-based approaches. 
Cross methods compute rank-$r$ approximation by interpolating the $n_1 \x n_2$ matrix on a cleverly chosen set of \emph{crosses}, proposed for instance in~\cite{tee-cross-2000,gt-maxvol-2001}. 
That requires $\O(nr)$ evaluations of matrix elements and can be used for matrices implicitly given by a subroutine that evaluates any prescribed element. 
However, the verification of the approximation requires several heuristics. 
If the matrix  is structured (sparse, Toeplitz or Hankel, low-rank, sum and/or product of above, etc.) and a fast matrix-by-vector product is available, then the Krylov-type methods are methods of choice with established convergence and complexity estimates. 

The generalization of the cross method to $3$-tensors (Cross3D, \cite{ost-tucker-2008}) requires the careful algorithmic implementation and several tweaks to make it really efficient. 
Generally, the Cross3D interpolates a tensor on $\O(nr+r^3)$ elements  and uses $\O(nr^2+r^4)$ additional operations. 
However, the accuracy check also requires heuristics. 

Recently, Eld\'en and Savas~\cite{eldensavas-krylov-2009} proposed two generalizations of Krylov methods to tensors: minimal and maximal Krylov recursion. 
The minimal recursion requires only $3r$ tensor-by-vector-by-vector multiplications to compute basis sets $U, V, W$ for the Tucker rank-$(r,r,r)$ approximation $(T),$ but sometimes it converges slowly and even does not converge (for example, for tensors with sufficiently different mode ranks).
The maximal recursion is convergent, but needs definitely unaffordable number of tensor operations. 
The goal of this paper is to propose algorithms that use tensor only through the tensor-by-vector-by-vector multiplication subroutine, have asymptotical complexity equal to the one of minimal Krylov recursion whereas have better convergence properties. 
We also try to keep a link with the existing theory of matrix approximation and with ideas used in~\cite{ost-tucker-2008} for the cross approximation of tensors.

This paper is organized as follows.
In the section~\ref{OPT} we propose the optimization of minimal Krylov recursion, inspired by the idea of maximization of the orthogonal component of new vectors.
In the section~\ref{MAT} we recall the Wedderburn rank reduction formula that gives a nice framework for many matrix factorizations~\cite{chu-rank1-1995}.
Then we propose a  variant of a matrix decomposition that is similar to the Gaussian elimination with pivoting but uses the matrix through the matrix-by-vector products.
In the section~\ref{TEN} we generalize this method to  the tensor case. 
The Wedderburn elimination process gives some freedom in a selection of the vectors to eliminate, and we use it proposing several `pivoting' strategies, that lead to different complexity estimates and convergence properties. 
Some pivoting strategies turn the proposed method into the minimal Krylov recursion and the optimized minimal Krylov recursion.
In the section~\ref{NUM} we apply the proposed family of algorithms for the approximation of different structured tensors and give numerical comparison with previous methods.
 
The \emph{tensor-by-vector-by-vector multiplication}, shortly the \emph{tenvec} operation, can be defined via the tensor-by-matrix products. 
For example, tenvec of $n_1\x n_2 \x n_3$ tensor $\A=[a_{ijk}]$ at modes $2,3$ reads
$$
u = \A \x_2 v^\t \x_3 w^\t, \qquad u_i = \Sum_{j=1}^{n_2} \Sum_{k=1}^{n_3} a_{ijk} v_j w_k
$$
where $v,w$ are vectors of length $n_2$ and $n_3$ respectively and the result is a vector $u$ of length $n_1.$
This is the solely tensor operation in algorithms in this paper and we propose extremely simple notation 
\begin{equation}\nonumber 
\A \x_2 v^\t \x_3 w^\t \eqdef \A v w, \qquad \A \x_3 w^\t \x_1 u^\t \eqdef \A w u, \qquad \A \x_1 u^\t \x_2 v^\t \eqdef \A u v,
\end{equation} 
relaxing the information on contraction modes, that is always obvious from the symbolic notation and mode sizes of vectors.
This simple notation (that generalizes the \emph{Tdd} notation~\cite{ali-krylov-1996}) reveals the analogy of tensor algorithms with the matrix case.

We consider the approximation of matrices and tensors in the Frobenius norm
$$
\|\A \|_F^2 \eqdef \f{\A}{\A}, \qquad \f{\A}{\B} \eqdef \Sum_{i=1}^{n_1} \Sum_{j=1}^{n_2} \Sum_{k=1}^{n_3} a_{ijk} b_{ijk}.
$$
For theoretical estimates we also use the spectral norm of tensor (cf.~\cite{defant-1993})
$$
\|\A\|_\tn \eqdef \max_{\|x\|=\|y\|=\|z\|=1} \A \x_1 x^\t \x_2 y^\t \x_3 z^\t = \max_{\|x\|=\|y\|=\|z\|=1} \f{\A}{x \o y \o z},
$$
induced by the standard vector norm $\|x\|^2 \eqdef \|x\|_2^2  =  (x,x) = \sum_{i=1}^{n} |x_i|^2.$ 

It is worth to mention that discussed algorithms first aim to construct the orthogonal bases $U, V$ and $W$ of the dominant mode subspaces, and the core $\G$ of the Tucker approximation can be computed afterwards.
For the fixed $U, V, W$ the most accurate approximation in Frobenius norm is obtained with the core 
\begin{equation}\label{eq:core}
 \G = \A \x_1 U^\t \x_2 V^\t \x_3 W^\t.
\end{equation}
In general,~\eqref{eq:core} can be computed using $r^2$ tenvecs, but faster implementations are available for certain structured tensors and can essentially improve the complexity. 
We can also be satisfied with the sub-optimal but fast formula for $\G,$ for example the one based on interpolation on a maximum-volume set of indices~\cite{gostz-maxvol-2010} as proposed in~\cite{ost-chem-2009}.
If the computation of~\eqref{eq:core} is not a necessary part of the discussed algorithm, this cost is not included in the total complexity.

\section{Minimal Krylov recursion and its optimization} \label{OPT}

The generalization of the Krylov subspace method to the problem of tensor approximation was first proposed by Eld\'en and Savas in~\cite{eldensavas-krylov-2009}. Two variants are discussed, namely the minimal and the maximum Krylov recursion.

The minimal Krylov recursion (MKR, see Alg.~\ref{alg:min}) requires $3$ tenvecs on each iteration, therefore, basis sets $U=U_r, V=V_r, W=W_r$ are computed in $3r$ tenvecs. 
However, their `quality', i.e. the accuracy of the approximation $\tilde\A = \G \x_1 U \x_2 V \x_3 W $ with the optimal core~\eqref{eq:core} may be low, and for some cases $\|\A - \tilde \A \|_F$ will not reduce at all, even for large ranks.

\begin{algorithm}[t]
\caption{\cite{eldensavas-krylov-2009} Minimal Krylov recursion for tensor approximation (MKR)}\label{alg:min}
\begin{algorithmic}
\REQUIRE Tenvec subroutine for tensor $\A$ 
\ENSURE Mode subspaces $U, V, W$ for the Tucker approximation
\item[\textbf{Initialization:}] Unit vectors $u_1, v_1$ 
\STATE  $w_1 := \A u_1 v_1/ \|\A u_1 v_1\|, \quad U_1 = [u_1], V_1=[v_1], W_1=[w_1]$ 
\FOR{$k=1,2,\ldots$}
     \STATE $u := \A v_k w_k;           \:\: u' := (I-U_k U_k^\t)u; \:\: u_{k+1} := u'/\|u'\|; \:\: U_{k+1} := [U_k \: u_{k+1}]$
     \STATE $v := \A w_k u_{k+1};       \:\: v' := (I-V_k V_k^\t)v; \:\: v_{k+1} := v'/\|v'\|; \:\: V_{k+1} := [V_k \: v_{k+1}]$
     \STATE $w := \A u_{k+1} v_{k+1};   \:\: w' := (I-W_k W_k^\t)w; \:\: w_{k+1} := w'/\|w'\|; \:\: W_{k+1} := [W_k \: w_{k+1}]$
\ENDFOR
\end{algorithmic}
\end{algorithm}

As an example, consider an $\nnn$ tensor $\A$ with only two non-zero slices%
\footnote{Here and after we use the MATLAB-style notation, with ``$:$'' denoting all possible index values}
$$
\A(:,:,1) = A_1, \qquad \A(:,:,2) = A_2.
$$ 
Obviously, if $A_1 \neq A_2,$ the mode-$3$ subspace of $\A$ consists  of two vectors $W = [e_1 \: e_2].$ 
Starting the MKR Alg.~\ref{alg:min} from some $u_1, v_1,$ we accumulate $W_2 = [e_1 e_2]$ in two steps. Then all the computed vectors $w$ have zero component $w'$ orthogonal to $W_2.$ 
This situation is referred to as \emph{breakdown}, since we cannot continue the process.

In~\cite{eldensavas-krylov-2009} Eld\'en and Savas propose to fix breakdowns by taking an orthogonal vector to the subspace $W.$ 
However with $w_3:=e_3$ we come to zero vector $u := \A v_3 w_3$ on the next iteration and can not continue the process. 
Another possible workaround, proposed in~\cite{eldensavas-krylov-2010} is to use the last basis vector on all the subsequent iterations, setting $w_3 := w_2, w_4 := w_2$ and so on. 
This should be done when the subspace of the mode-3 vectors of $\A$ equals to $\Span W_k.$ 
However, it is difficult to check this fact using only tenvec operations and it is not clear how to do this in the non-exact case which is always true in the machine precision arithmetic.
In the numerical examples (see section~\ref{NUMc}) we show that the poor convergence and even the stagnation can occur also for `more practical' situations.

From this discussion we derive the idea of optimization of the Alg.~\ref{alg:min}. Consider the step for $u$ and generalize it as follows.
$$
u := \A (V_k \vv_k) (W_k \ww_k); \quad u' := (I-U_k U_k^\t)u; \quad u_{k+1} := u'/\|u'\|; \quad U_{k+1} := [U_k \: u_{k+1}].
$$
New direction $u$ is generated by the tenvec of $\A$ with \emph{some} vector from $\Span V$ and \emph{some} vector from $\Span W.$
In MKR we always take $\vv_k = \ww_k = e_k,$ but there is no background theory behind this choice and it seems to be not optimal in practice.
To improve the approximation properties of $U_k,$ we could choose $\vv_k$ and $\ww_k$ to maximize the norm of the orthogonal component $u'.$
Since 
$$
u' = (I - U_k U_k^\t) u = \left( \A \x_1 (I - U_k U_k^\t) \x_2 V_k^\t \x_3 W_k^\t \right) \x_2 \vv_k^\t \x_3 \ww_k^\t,
$$
we are to find 
\begin{equation}\label{eq:optmax}
\vv_k, \ww_k = \arg\max_{\|\vv\|=1, \|\ww\|=1}  \| \B \vv \ww \| 
\qquad \mbox{for} \quad 
\B=\A \x_1 (I - U_k U_k^\t) \x_2 V_k^\t \x_3 W_k^\t.
\end{equation}

\begin{algorithm}[t]
\caption{\cite{tuckerals-1980}~ALS rank-$(1,1,1)$ iteration}\label{alg:als} 
\begin{algorithmic}
\REQUIRE Tenvec subroutine for tensor $\A$ 
\ENSURE Best rank-$(1,1,1)$ approximation $\sigma u \o v \o w$ for tensor
\item[\textbf{Initialization:}] Unit vectors $v, w$
\FOR{$k=1, \ldots, \pals$}
      \STATE  $u:=\A v w; \quad \sigma := \|u\|; \quad u := u/\|u\|$
      \STATE  $v:=\A w u; \quad \sigma := \|v\|; \quad v := v/\|v\|$
      \STATE  $w:=\A u v; \quad \sigma := \|w\|; \quad w := w/\|w\|$
\ENDFOR
\end{algorithmic}
\end{algorithm}

The global maximum is not required. We can satisfy with the sufficiently large $u'$ that can be found by several iterations of the \emph{alternating least squares} method (ALS, see~\cite{tuckerals-1980} and Alg.~\ref{alg:als}), for which the local linear convergence is proved in the rank-$(1,1,1)$ case~\cite{lathauwer-rank1-2000,zhang-rank1-2001}. 
If the ALS converges to the best rank-one approximation, it also solves~\eqref{eq:optmax} as a dual problem.
To see this, consider the approximation $\tilde\B = b u \o \vv \o \ww$ with unit $u, \vv, \ww$ and compute the optimal $1\x1\x1$ core $b$ by~\eqref{eq:core}
$$
b:=\B \x_1 u^\t \x_2 \vv^\t \x_3 \ww^\t = \f{\B}{u\o\vv\o\ww} = (\B\vv\ww, u).
$$
The approximation $\tilde\B$ that minimizes the approximation error
$$
\|\B - \tilde \B\|_F^2 = \|\B\|_F^2 - 2 \f{\B}{b u\o\vv\o\ww} + \|b u\o\vv\o\ww\|^2_F = \|\B\|^2_F - |b|^2
$$
also maximizes $|b|$ and solves~\eqref{eq:optmax}, since 
$$
\max_{\|u\|=\|\vv\|=\|\ww\|=1} |b| = 
\max_{\|\vv\|=\|\ww\|=1} \max_{\|u\|=1} (\B\vv\ww, u) = \max_{\|\vv\|=\|\ww\|=1} \|\B\vv\ww\|.
$$

Each ALS iteration requires $3$ tenvecs with $\B,$ that can be rendered by $3$ tenvecs with $\A$ and $\O(nk)$ operations for the orthogonalization.
The optimization of the MKR is summarized in terms of tenvecs operations in the Alg.~\ref{alg:opt}. 

\begin{algorithm}[t]
\caption{Optimized minimal Krylov recursion for tensor approximation}\label{alg:opt}
\begin{algorithmic}[1]
\REQUIRE Tenvec subroutine for tensor $\A,$  tolerance parameter $\tol$
\ENSURE Mode subspaces $U, V, W$ for approximation $\A \approx \tilde \A = \G \x_1 U \x_2 V \x_3 W$
\item[\textbf{Initialization:}] Unit vectors $u_1, v_1$ 
\STATE  $w_1 := \A u_1 v_1/\|\A u_1 v_1\|, \quad U_1 = [u_1], V_1 = [v_1], W_1 = [w_1]$ 
\STATE $\updu=\updv=\updw=\TRUE; \quad k=l=m=1$ 
\WHILE{$\updu$ or $\updv$ or  $\updw$}
      \IF[Proceed with the new vector $u$ if required]{$\updu$}
      \STATE            Define $\B=\A \x_1 (I - U_k U_k^\t) \x_2 V_l^\t \x_3 W_m^\t$ 
      \STATE            Find $\vv_k, \ww_k := \arg\max_{\|\vv\|=1, \|\ww\|=1}\|\B \vv \ww\|$ by $\pals$ ALS steps, see Alg.~\ref{alg:als}
      \STATE            $u := \A (V_l\vv_k) (W_m\ww_k); \quad u' := (I-U_k U_k^\t)u$
                  \IF[Breakdown]{$\|u'\| < \tol \|u\|$}
                      \STATE   Fix breakdown or set $\updu:=\FALSE$ 
                  \ELSE
                      \STATE          $u_{k+1}:=u'/\|u'\|; \quad U_{k+1} := [U_k \: u_{k+1}]; \quad k := k+1$ 
                  \ENDIF
      \ENDIF
      \STATE \COMMENT{Proceed with new vector $v$ if required}
      \STATE \COMMENT{Proceed with new vector $w$ if required}
\ENDWHILE
\end{algorithmic}
\end{algorithm}

If $\pals$  ALS iterations are used to solve~\eqref{eq:optmax} on Step $6,$ then Alg.~\ref{alg:opt} requires $(3+9\pals)r$  tenvecs and $\O(nr^2\pals)$ additional operations.

In numerical experiments (see section~\ref{NUM}) we show that the optimized MKR Alg.~\ref{alg:opt} shows the better convergence that the MKR Alg.~\ref{alg:min}.
However, convergence estimates are still missing, even in the exact low-rank case. 
Is it possible to develop a method that has guaranteed convergence at least in the exact case? 
We start from the matrix case, using Wedderburn rank-one reduction formula.
Then, in section~\ref{TEN}, we propose generalization of Wedderburn method to tensor case and reintroduce the MKR and the optimized MKR as different versions of the same Wedderburn process, for which we discuss the convergence and give numerical examples and comparison.

\section{Matrix approximation using Wedderburn rank reduction formula} \label{MAT}
\subsection{Preliminaries}
Many matrix decomposition algorithms can be represented as a sequence of rank-one \emph{Wedderburn updates}~\cite{wedderburn}.
For a matrix $A$ and vectors $x,y$ of appropriate sizes, such that $x^\t A y \neq 0,$ the matrix
\begin{equation}\label{eq:wdm0}
B = A  - \frac{A y x^\t A}{x^\t A y}.
\end{equation}
has $\rank B = \rank A - 1.$
For the rank-$r$ matrix $A_0=A$ after $r$ updates of form
\begin{equation}\label{eq:wdm}
A_k = A_{k-1} - \frac{A_{k-1} y_k x_k^\t A_{k-1}}{x_k^\t A_{k-1} y_k}
\end{equation}
with $\omega_k \eqdef x_k^\t A_{k-1} y_k \neq 0,$  the matrix $A_r$ becomes zero and the rank-$r$ decomposition of $A$ can be constructed.

In \cite{chu-rank1-1995} the properties of the Wedderburn sequence~\eqref{eq:wdm} are studied in much detail. 
We recall a list of basic facts about~\eqref{eq:wdm} in the compact and `more matrix' form. 

\begin{statement}\label{stat1}
Matrix $A_k$ writes as $A_k = P_k^\t A = A Q_k$ with $P_0 = I,$ $Q_0=I$ and
\begin{equation}\label{eq:wdm1}
P_k = P_{k-1} - \omega_k^{-1} P_{k-1} x_k y_k^\t A^\t P_{k-1}, \qquad
Q_k = Q_{k-1} - \omega_k^{-1} Q_{k-1} y_k x_k^\t A    Q_{k-1}.
\end{equation}
\end{statement}

\begin{corollary}\label{cor1}
$P_k x_k = 0$ and  $Q_k y_k  = 0.$
\end{corollary}

\begin{statement}\label{stat2}
Matrices $P_k, Q_k$ are projectors, i.e. $P_k^2 = P_k$ and $Q_k^2 = Q_k.$
\end{statement}
Obviously, $P_0 = P_0^2 = I,$ $Q_0 = Q_0^2 = I$ and from $P_{k-1}^2=P_{k-1}$ it follows that
\begin{equation}\nonumber
\begin{split}
P_k^2  = P_{k-1}^2 - \omega_k^{-1} P_{k-1}^2 y_k x_k^\t A^\t P_{k-1} & - \omega_k^{-1} P_{k-1} x_k y_k^\t A^\t P_{k-1}^2 + {} \\ &
                     + \omega_k^{-2} P_{k-1} x_k  \underbrace{y_k^\t A^\t P_{k-1}^2 x_k}_{\omega_k} y_k^\t A^\t P_{k-1} = P_k.
\end{split}
\end{equation}

\begin{statement}\label{stat3}
$P_k  P_m = P_m P_k = P_{\max(m,k)},$ and  $Q_k  Q_m = Q_m Q_k = Q_{\max(m,k)}.$ 
\end{statement}
Start from 
$$
Q_k Q_{k-1} = Q_{k-1}^2 - \omega_k^{-1}Q_{k-1}y_kx_k^\t A Q_{k-1}^2 = Q_{k-1} - \omega_k^{-1}Q_{k-1}y_kx_k^\t A Q_{k-1} = Q_k
$$ 
and complete the proof by induction.

\begin{corollary}\label{cor2}
For $m \leq k$ it holds $P_k x_m = P_k P_m x_m = 0$ and  $Q_k y_m = Q_k Q_m y_m = 0.$
\end{corollary}

\begin{statement}\label{stat4}
For \emph{biconjugate vectors} $u_k \eqdef P_{k-1} x_k$ and $v_k \eqdef Q_{k-1} y_k$ it holds
$$
u_k^\t A v_k = x_k^\t A_{k-1} y_k \eqdef \omega_k,
$$
since $A_{k-1} = A Q_{k-1} =  A Q_{k-1}^2  = A _{k-1} Q_{k-1} = P_{k-1}^\t A Q_{k-1}.$
Also for $m \neq k$ it holds
$$
u_m^\t A v_k = x_m^\t P_{m-1}^\t A Q_{k-1} y_k =  x_m^\t A_{\max (m-1,k-1) } y_k = 0.
$$
With $U_k \eqdef [u_1 \:\ldots\: u_k],$  $V_k \eqdef [v_1 \:\ldots\: v_k]$ and $\Omega_k \eqdef \diag(\omega_1 \:\ldots\: \omega_k)$ we conclude
\begin{equation}\label{eq:wdm6}
U^\t_k A V_k = \Omega_k.
\end{equation}
\end{statement}

\begin{corollary}\label{cor3}
For the valid Wedderburn process $U_k$ and $V_k$ have full rank.
\end{corollary}
\begin{corollary}\label{cor4}
For $m \leq k$ it holds $P_k u_m = P_k P_{m-1} x_m = P_k x_m = 0$ and  $Q_k v_m = 0.$
\end{corollary}

\begin{statement}\label{stat5}
Each rank elimination step~\eqref{eq:wdm} writes without multiplication by $A_{k-1}$ as follows
$$
A_k = A_{k-1} - \frac{A v_k u_k^\t A}{u_k^\t A v_k} = A - \Sum_{p=1}^k A v_p \omega_p^{-1} u_p^\t A.
$$
Now the rank-$k$ \emph{Wedderburn approximation}  reads
\begin{equation}\label{eq:wda}
 \tilde A_k \eqdef A V_k \Omega_k^{-1} U_k^\t A,  \qquad
 A = A_k + \tilde A_k.
\end{equation}
For the rank-$r$ matrix $A$ the residual $A_r = 0$ and the approximation $\tilde A_r$ is exact.
\end{statement}

We see that each Wedderburn update~\eqref{eq:wdm} adds vector $y_k$ to the kernel and $x_k$  to the cokernel of the residual $A_k = A - \tilde A_k.$ 
The approximation $\tilde A_k,$ as a linear operator, interpolates $A$ on subspaces spanned by $X_k$ and $Y_k,$ exactly
\begin{equation}\label{eq:wdm9}
 X_k^\t \tilde A_k = X_k^\t A, \quad \tilde A_k Y_k = A Y_k, \qquad 
 U_k^\t \tilde A_k = U_k^\t A, \quad \tilde A_k V_k = A V_k.
\end{equation}
This  can be associated with the \emph{Gaussian elimination}, that gives an approximation exact on certain rows and columns of the matrix. 
In this respect we refer to the process~\eqref{eq:wdm} as to the \emph{Wedderburn elimination.}

In \cite{chu-rank1-1995} it is shown how the proper choice of $x_k, y_k$ reduces the Wedderburn elimination to the well-known matrix decompositions such as  LU, QR, cross factorizations and Lanczos bidiagonalization.
We consider another principle for the selection of vectors $x_k, y_k$ at each step in the matrix case, that can be associated with Gaussian elimination with column or row pivoting. It produces a new method for matrix approximation, that is simply generalized to tensors while maintaining the convergence.

\subsection{Pivoting in Wedderburn elimination}
The idea behind the proposed choice of $x_k, y_k$ is minimization of Frobenius norm of the residual, that is important  when we deal with the full-rank matrix of a low $\eps$-rank, i.e. that can be approximated by the low-rank matrix with the accuracy $\eps.$
As shown in~\cite{chu-rank1-1995}, the minimization of residual w.r.t. unit $x_k, y_k$ gives $k$-th singular vectors of $A,$ that are not known in advance. We propose another minimization strategy, based on the following theorem.

\begin{theorem}\label{thm1}
Consider Wedderburn step~\eqref{eq:wdm0}. Then
\begin{equation}\label{eq:optx}
\mbox{for fixed $y,$} \qquad 
      x_{\mathrm{opt}} \eqdef \arg \min_{\|x\|=1} \left\| A - \frac{A y x^\t A}{x^\t A y} \right\|_F = \frac{ A y }{\|Ay\|}; 
\end{equation}
\begin{equation}\label{eq:opty}
\mbox{for fixed $x,$} \qquad 
      y_{\mathrm{opt}} \eqdef \arg \min_{\|y\|=1} \left\| A - \frac{A y x^\t A}{x^\t A y} \right\|_F = \frac{ A^\t x}{\|A^\t x\|}.
\end{equation}
\end{theorem}
\begin{proof}
Since the valid Wedderburn step does not depend on scaling of $x, y,$ for the fixed $y$ we can constrain $x$ to satisfy $x^\t A y = 1.$ Then
$$
x_{\mathrm{opt}} = \min_{x:\: x^\t A y = 1 } \left\|(A y)  (A^\t x) - A \right\|_F^2 =
                   \min_{x:\: x^\t A y = 1 } \left\|(A^\t x) (A y)^\t - A^\t\right\|_F^2.
$$
The least squares problem solves by $(A^\t x)_\mathrm{opt} = A^\t A y / \|Ay\|^2$  and $x_\mathrm{opt} = Ay / \|Ay\|^2.$ 
After normalization we have~\eqref{eq:optx}, and~\eqref{eq:opty} follows by substituting $A = A^\t.$
\end{proof}

Eqs.~\eqref{eq:optx} and~\eqref{eq:opty} show how to reach the fast decay of the residual in the Wedderburn process~\eqref{eq:wdm} either by choosing optimal $x_k$ for given $y_k$ or by choosing optimal $y_k$ for given $x_k.$ This can be associated with the column and row pivoting in the Gaussian elimination. 
Thus, we refer to  the Wedderburn process~\eqref{eq:wdm} with arbitrary $y_k$ and optimal $x_k$ as to \emph{Wedderburn elimination with column pivoting} (WCP)  and to 
Wedderburn process with arbitrary $x_k$ and optimal $y_k$ as to \emph{Wedderburn elimination with row pivoting} (WRP).
The steps of WCP and WRP shortly read
$$
\mbox{choose $y_k,$ set} \quad x_k = \frac{A_{k-1} y_k}{ \|A_{k-1} y_k\| }, \quad 
A_k  =  \left( I - x_k x_k^\t \right) A_{k-1};
\eqno{(WCP)}
$$
$$
\mbox{choose $x_k,$ set} \quad y_k = \frac{A_{k-1}^\t x_k}{ \|A_{k-1}^\t x_k\| }, \quad 
A_k  =  A_{k-1} \left( I - y_k y_k^\t \right). 
\eqno{(WRP)}
$$

\begin{theorem}\label{thm2}
For the Wedderburn elimination with column pivoting it holds 
\begin{enumerate}
      \item $X_k = [x_1 \ldots x_k]$ has orthonormal columns $X_k^\t X_k = I;$ 
      \item $P_k=P_k^\t = I - X_k X_k^\t$ is the projector on the subspace orthogonal to $\Span X_k;$
      \item biconjugate vectors $u_k \eqdef P_{k-1}x_k = x_k.$
\end{enumerate}
For the Wedderburn elimination with row pivoting  it holds 
\begin{enumerate}
      \item$ Y_k = [y_1 \ldots y_k]$ has orthonormal columns $Y_k^\t Y_k = I;$ 
      \item $Q_k=Q_k^\t = I - Y_k Y_k^\t;$
      \item $v_k \eqdef Q_{k-1} y_k = y_k.$
\end{enumerate}
\end{theorem}
\begin{proof}
Let us prove statements for WRP by induction. It is easy to check them for $k=1.$ 
Suppose they hold at step $k-1.$ Optimal choice of $y_k$ by~\eqref{eq:opty} reads
$$
y = A_{k-1}^\t x_k = Q_{k-1}^\t A^\t x_k = (I - Y_{k-1} Y_{k-1}^\t) A^\t x_k, \qquad y_k = {y}/{\|y\|}.
$$ 
This shows $\|y_k\|=1,$  $Y_{k-1}^\t y_k = 0$ and the first statement follows for $Y_k := [Y_{k-1} \: y_k].$ 
By substituting $y_k=A_{k-1}^\t x_k/\|A_{k-1}^\t x_k\|$ in~\eqref{eq:wdm1} we prove the second statement
\begin{equation}\nonumber
\begin{split}
Q_k & = Q_{k-1} - \frac{ Q_{k-1} y_k x_k^\t A_{k-1} }{ x_k^\t A_{k-1} y_k } = Q_{k-1} \left(I - \frac{y_k y_k^\t}{y_k^\t y_k} \right) = \\
    & = \left(I - Y_{k-1} Y_{k-1}^\t \right) \: \left( I - y_k y_k^\t \right) = I - Y_{k-1} Y_{k-1}^\t - y_k y_k^\t = I - Y_k Y_k^\t.
\end{split}
\end{equation}
Finally, the last statement reads
$$
v_k \eqdef Q_{k-1} y_k = Q_{k-1} \frac{ Q_{k-1}^\t A^\t x_k }{\|A_{k-1}^\t x_k\|}  = 
      \frac{ Q_{k-1}^2 A^\t x_k }{\|A_{k-1}^\t x_k\|} = \frac{ Q_{k-1}^\t A^\t x_k }{\|A_{k-1}^\t x_k\|} = y_k.
$$
Statements for the WCP are proven in the same way.
\end{proof}

\begin{algorithm}[t]
\caption{Wedderburn elimination with column pivoting (WCP)} \label{alg:wcp}
\begin{algorithmic}[1]
\REQUIRE Matvec subroutine for the matrix $A,$  tolerance parameter $\tol,$ accuracy $\eps$
\ENSURE Approximation $\tilde A$ with the accuracy $\|A - \tilde A\|_F \lesssim \eps\|\tilde A\|_F$
\item[\textbf{Initialization:}] $k=0,$ $X_0 = [\varnothing],$ $B_0 = [\varnothing],$ $\nrm=0$ 
\REPEAT 
      \STATE $k:=k+1,$ choose unit vector $y_k$ 
      \STATE $x  :=  A y_k; \quad x' := (I - X_{k-1} X_{k-1}^\t) x$
      \IF[Breakdown]{$\|x'\| < \tol \|x\|$}
            \RETURN $\tilde A = X_{k-1} B_{k-1}^\t$ \COMMENT{or repeat current iteration with another $y_k$}
      \ELSE      
            \STATE $x_k := x' / \|x'\|$
      \ENDIF
      \STATE $b_k := A^\t x_k, \quad \err := \|b_k\|, \quad \nrm^2 := \nrm^2 + \|b_k\|^2$ 
      \STATE $X_k := [X_{k-1} \: x_k], \quad B_k := [B_{k-1} \: b_k]$
\UNTIL{$\err \leq \eps \nrm$}
\RETURN $\tilde A = X_k B_k^\t$
\end{algorithmic}
\end{algorithm}

Based on this theorem, we propose a Wedderburn elimination algorithm with the column pivoting (WCP, see Alg.~\ref{alg:wcp}).
The approximation is sought in the form $\tilde A_k = X_k B_k^\t$ with $B_k = A^\t X_k,$ that follows from~\eqref{eq:wdm9} and orthogonality of biconjugate vectors $U_k = X_k,$ result of the Theorem~\ref{thm2}. 
To explain the stopping criteria of Alg.~\ref{alg:wcp}, write
$$
\nrm \eqdef \|\tilde A_k\|_F = \| X_k X_k^\t A \|_F = \|A^\t X_k \|_F = \|B_k\|_F, \quad
\err \eqdef \|\tilde A_k - \tilde A_{k-1}\|_F = \|b_k\|_F,
$$
where the first part estimates the norm of matrix by the norm of approximation and the second estimates the error of approximation by the norm of the update.
The \emph{breakdown} can occur if new vector $x$ has the neglectable component $x'$ orthogonal to the accumulated subspace~$X_{k-1}.$
This is regulated by the tolerance parameter $\tol,$ that could be chosen close to machine precision. 
We can try to fix the  breakdown by another selection of $y_k,$ or choose to terminate the algorithm.

The WRP version of algorithm follows by substituting $A=A^\t.$ 

\subsection{Relation to SVD and Lanczos bidiagonalization}
On each step of the WCP Alg.~\ref{alg:wcp} we perform the Wedderburn elimination, choosing optimal `pivot' $x_k$ for given $y_k.$ 
It is also important to select proper `leading vectors' $y_k$ to obtain faster convergence of approximation and avoid breakdowns.
We could think about the maximization of $\omega_k = x_k^\t A_{k-1} y_k = \|A_{k-1} y_k\| = \|x'\|.$ Solving exactly
\begin{equation}\label{eq:best}
 y_k = \arg \max_{\|y\|=1} \|A_{k-1} y\| = \arg \max_{\|y\|=1} \left\| (I - X_{k-1} X_{k-1}^\t) A y \right\|,
\end{equation}
we reduce the Wedderburn process to a sequence of best rank-one approximations with $\omega_k$ and  $x_k, y_k$  being singular values of $A$ (sorted descending) and corresponding left and right singular vectors.
Therefore, the WCP with selection of the leading vector by~\eqref{eq:best} is equivalent to the SVD and provides the best rank-$r$ approximation in the Frobenius norm.
This approach can be associated with the \emph{full pivoting} in the Gaussian elimination.
Each maximization problem~\eqref{eq:best} can be accurately solved by power iterations using only matvec operations, but
in matrix case this approach generally is considered as quite expensive and faster alternatives are used.

One of them in the Lanczos bidiagonalization, see~\cite{golubkahan-1965} and Alg.~\ref{alg:lanczos}.
It generates bases $X_k = [x_1 \ldots x_k]$ and $Y_k = [y_1 \ldots y_k]$ such that $X_k^\t X_k =  I,$ $Y_k^\t Y_k = I$ and the matrix $X_k^\t A Y_k$ is bidiagonal.

The following theorem shows that the Lanczos bidiagonalization is similar to the Alg.~\ref{alg:wcp} if the leading vector is selected as follows
\begin{equation}\label{eq:choice}
y_{k+1} = \frac{A_{k-1}^\t x_k}{\|A_{k-1}^\t x_k\|} = \frac{A^\t x_k}{\|A^\t x_k\|}.
\end{equation}
This value is well defined, since
$$
A_{k-1}^\t x_k = A^\t P_{k-1} x_k = A^\t (I - X_{k-1} X_{k-1}^\t) x_k = A^\t x_k = b_k,
$$
and $\|A_{k-1}^\t x_k\| \eqdef \err$ vanishes only when the stopping criteria $\err < \eps \nrm$ is met and the next leading vector $y_{k+1}$ is not required.

\begin{algorithm}[t]
\caption{\cite{golubkahan-1965} Lanczos bidiagonalization} \label{alg:lanczos}
\begin{algorithmic}
\item[\textbf{Initialization:}]  $y_0=0, \beta_0=0,$  unit vector $x_0$ 
\FOR{$k=1,2,\ldots$}
      \STATE $y:= A^\t x_{k-1}; \quad y' := y - \beta_{k-1} y_{k-1}, \quad \alpha_k := \|y'\|,   \quad y_k := y' / \|y'\|$
      \STATE $x:= A y_k;        \quad x' := x - \alpha_k x_k,        \quad \beta_k := \|x'\|,    \quad x_k := x' / \|x'\|$
\ENDFOR
\end{algorithmic}
\end{algorithm}

\begin{theorem}\label{thm3}
Vectors $\Xgk = [\xgk_1,\ldots,\xgk_k]$ of the Lanczos process (Alg.~\ref{alg:lanczos}) initialized by $x_0,$ coincide with
vectors $\Xwcp_k = [\xwcp_1, \ldots, \xwcp_k]$ generated by WCP Alg.~\ref{alg:wcp} that starts from $y_1 = A^\t x_0$ and chooses $y_k$ as proposed by~\eqref{eq:choice},  providing both algorithms do not meet breakdowns. 
\end{theorem}
\begin{proof}
If no breakdowns are met, then
\begin{equation}\nonumber
      \begin{split}
      \Span \Xwcp_k &=   \Span\{ A y_1, (AA^\t) A y_1, \ldots, (AA^\t)^{k-1} A y_1 \}, \\
      \Span \Xgk_k  &=   \Span\{ (AA^\t) x_0, (AA^\t)^2 x_0, \ldots, (AA^\t)^k x_0 \}.
      \end{split}
\end{equation}
If $y_1 = A^\t x_0$ then $\Span \Xwcp_k = \Span \Xgk_k$ for every $k.$ 
Since columns of $\Xwcp_k$ are orthonormal (Theorem~\ref{thm2}), as well as $\Xgk_k,$ we conclude $\Xwcp_k = \Xgk_k.$ 
\end{proof}

Thus, in terms of $x_k,$  the WCP gives the same result as the Lanczos bidiagonalization. 
However, the sequence $y_k$ in the WCP differs from the one of the Lanczos bidiagonalization. First, we note that
$$
y = A_{k-1}^\t A_{k-1} y_k = A^\t P_{k-1} P_{k-1}^\t A y_k = A^\t P_{k-1} A y_k = Q_{k-1}^\t A^\t A y_k, \quad
y_{k+1} = \frac{y}{\|y\|}
$$
and $\Span Y_k = \Span\{ y_1, (A^\t A)y_1, \ldots, (A^\t A)^{k-1} y_1 \}.$
Sequence $y_k$ is `almost orthogonal', i.e. for $m \leq k-2$ it holds
$$
y_k^\t y_m = \frac{ y_{k-1}^\t A^\t A Q_{k-2} y_m }{\|A^\t A y_{k-1}\|} = 0 
$$
since $Q_{k-2} y_m = 0$ by the Corollary~\ref{cor2}.
Also, for $X_k$  and $Y_k$ generated by the WCP, the matrix $X_k^\t A Y_k$ is tridiagonal, i.e.
$$
x_m^\t A y_k  = y_{m+1}^\t y_k = 0, \qquad \mbox{for} \quad m \notin \{k-2, k-1, k\}.
$$
\begin{remark}
For $k \geq 3$ the vector $X_{k-1}^\t x = X_{k-1}^\t A y_k$  on Step~3 of Alg.~\ref{alg:wcp} has only two last nonzero components.
Hence, the orthogonalization step can be simplified to
$$
x : = A y_k, \quad x' : = (I - x_{k-2}x_{k-2}^\t- x_{k-1}x_{k-1}^\t)x.
$$
\end{remark}
This allows the short recursion for the orthogonalization of $x_k$ in the WCP Alg.~\ref{alg:wcp}, as it is done in the Lanczos bidiagonalization Alg.~\ref{alg:lanczos}.  
We summarize this version of the WCP in the Alg.~\ref{alg:wcp1}. 
In the machine arithmetic the orthogonality can be violated by roundoff errors and the reorthogonalization is required.
The convergence of the Lanczos bidiagonalization method is well studied and although this method can have breakdowns, it is considered to converge from `almost every' initial vector~\cite{paige-1976,golub-1996,greenbaum-1997}.

\begin{algorithm}[t]
\caption{WCP with Lanczos-like selection of leading vector} \label{alg:wcp1}
\begin{algorithmic}[1]
\REQUIRE Matvec subroutine for matrix $A,$  tolerance parameter $\tol,$ accuracy $\eps$
\ENSURE Approximation $\tilde A$ with accuracy $\|A - \tilde A\|_F \lesssim \eps\|\tilde A\|_F$
\item[\textbf{Initialization:}] $k=0, \quad x_0 = x_{-1} = 0, \quad \tilde A_0 = 0, \quad \nrm=0,$ unit vector $y_1$ 
\REPEAT 
      \STATE $k:=k+1, \quad x  :=  A y_k; \quad x' := (I - x_{k-2}x_{k-2}^\t-x_{k-1}x_{k-1}^\t) x$ 
      \IF[Breakdown]{$\|x'\| < \tol \|x\|$}
            \RETURN $\tilde A = \tilde A_{k-1}$ 
      \ELSE      
            \STATE $x_k := x' / \|x'\|$
      \ENDIF
      \STATE $y_{k+1} := A^\t x_k, \quad \err := \|y_{k+1}\|, \quad \nrm^2 := \nrm^2 + \|y_{k+1}\|^2$ 
      \STATE $\tilde A_k = \tilde A_{k-1} + x_k y_{k+1}^\t$
\UNTIL{$\err \leq \eps \nrm$}
\RETURN $\tilde A = \tilde A_k$
\end{algorithmic}
\end{algorithm}

\section{Tensor approximation using Wedderburn rank reduction} \label{TEN}

\subsection{Computing dominant mode subspaces by WCP algorithm}
We are ready to propose the extension of the Wedderburn elimination to the tensor case. 
For the Tucker approximation we are to approximate the dominant subspaces by $U, V, W$ that contain much of information about mode vectors of tensor.
In~\cite{lathauwer-svd-2000} this is done by SVD applied to the unfoldings of an $n_1 \x n_2 \x n_3$ tensor $\A=[a_{ijk}].$ 
They are matrices
\begin{equation}\label{eq:unf}
A^{(1)} = [a_{i, jk}^{(1)}], \quad
A^{(2)} = [a_{j, ki}^{(2)}], \quad
A^{(3)} = [a_{k, ij}^{(3)}], \qquad  a_{i, jk}^{(1)}=a_{j, ki}^{(2)}=a_{k, ij}^{(3)}=a_{ijk},
\end{equation}
of size $n_1 \x n_2n_3,$ $n_2 \x n_1n_3 $ and $n_3 \x n_1n_2,$ that consist of columns, rows and tube fibres of $\A,$ respectively. 
Left singular vectors after appropriate truncation give Tucker factors $U, V, W,$ and the core is found by~\eqref{eq:core}. 
Since the SVD is applied to tensors that are given as full array of elements, the computation costs $\O(n^4)$ and  can not is not possible for large tensors (even sparse or structured).  
We approximate dominant subspaces of mode vectors applying the Wedderburn elimination to the unfoldings.

We aim for algorithms that compute the rank-$(r_1, r_2, r_3)$ approximation of an $n_1 \x n_2 \x n_3$ tensor using $\O(r^2)$ tenvecs and $\O(n)$ additional operations. This is asymptotically equal to the cost of the MKR Alg.~\ref{alg:min}, if the core tensor is computed.
With this restriction for $n_1 \x n_2n_3$ unfolding $A = A^{(1)} = [a_{i;\: jk}] $ we can not compute $x = A \mathbf{y}$ and $\mathbf{y}=A^\t x$ for arbitrary vectors $\mathbf{y}$ of size $n_2n_3$ and $x$ of size $n_1,$ since these operations require $\O(n^2)$ storage and $n$ tenvecs.
To develop the Krylov-type methods and stay within the linear complexity in mode size, we should use only those operations, that can be accomplished by small number of tenvecs. For instance, $x = A\mathbf{y}$ is substituted by 
$$
x = A (y \o z) = \A \x_2 y^\t \x_3 z^\t = \A yz, 
$$
that means that we will use only those `long vectors' $\mathbf{y}$ that are the tensor product $y \o z$ of some $y$ of size $n_2$ and  $z$ of size $n_3.$ Evaluation of $\mathbf{y} = A^\t x = \A \x_1 x$ is also  infeasible and we substitute it with the approximation
$$
y \o z \approx \A \x_1 x,
$$
that means that we develop the algorithms that give accurate tensor approximation, using only certain rank-one approximation of $\A\x_1 x$ instead of the precise result.

In the matrix case the difference between column and row pivoting in the Wedderburn elimination is not significant, since the properties of column basis in the WCP coincide with ones of row basis in the WRP and vice versa.
In the tensor case with imposed restriction the multiplication $A(y \o z) $ is accurate, but the multiplication $A^\t x = \A\x_1 x^\t$ is always done approximately with an error of truncation to rank one. 
Considering this difference, we use the WCP algorithm for the tensor case, since good properties of $X_k$ given by Theorem~\ref{thm2} persist here.
Using the direct analogy with the WCP Alg.~\ref{alg:wcp} for matrices, we propose a method to derive dominant mode-$1$ subspace $U,$ see Alg.~\ref{alg:wt1}.

\begin{algorithm}[t]
\caption{Wedderburn elimination for dominant subspace computation} \label{alg:wt1}
\begin{algorithmic}[1]
\REQUIRE Tenvec subroutine for $\A,$ tolerance $\tol, $ accuracy $\eps,$ maximum size $\rmax.$
\ENSURE Mode subspace $U$ for Tucker approximation $\tilde\A$ such that  $\|\A - \tilde\A\|_F \lesssim \eps\|\tilde\A\|_F$
\item[\textbf{Initialization:}] $X_0 = [\varnothing], \quad \nrm=0, \quad k=0$ 
\REPEAT
      \STATE $k:=k+1,$ choose unit vectors $y_k, z_k$ (see Section~\ref{LEAD})
      \STATE $x  :=  \A y_k z_k; \quad x' := (I - X_{k-1} X_{k-1}^\t) x$
      \IF[Breakdown]{$\|x'\| < \tol \|x\|$}
            \RETURN $U = X_{k-1}$ \COMMENT{or repeat current iteration with another $y_k, z_k$}
      \ELSE      
            \STATE $x_k := x' / \|x'\|; \quad X_k := [X_{k-1} \: x_k]$
      \ENDIF
      \STATE Choose unit $z$ randomly
      \FOR[Power iterations to approximate $\A \x_1 x_k^\t \approx: \sigma y z^\t $]{$k=1,\ldots,\ppow$}
       \STATE $y := (\A \x_1 x_k^\t)^{\phantom{\t}} z    = \A \x_1 x_k^\t \x_3 z^\t = \A z x_k, \quad \sigma := \|y\|, \quad y := y/\|y\|$
       \STATE $z := (\A \x_1 x_k^\t)^\t y = \A \x_1 x_k^\t \x_2 y^\t = \A x_k y,                \quad \sigma := \|z\|, \quad z := z/\|z\|$
      \ENDFOR
      \STATE $\err:=\sigma, \quad \nrm^2 := \nrm^2 + \err^2$
\UNTIL{$\err \leq \eps \nrm$ or $k = \rmax$}
\RETURN $U = X_k$
\end{algorithmic}
\end{algorithm}

Each iteration begins with the choice of the `long' leading vector in the form $y_k \o z_k.$ It can be done arbitrarily, and we propose some good strategies in the section~\ref{LEAD}. The new direction $x$ and basis vector $x_k$ appear exactly like in the matrix case.
To terminate the method, we can use one or both of the stopping criteria: 
\begin{itemize}
 \item fix maximum number of iterations, i.e. desired size of basis $U$ by $\rmax,$ 
 \item find basis $U$ that allows approximation of $\A$ with relative accuracy $\eps.$
\end{itemize} 
In the WCP Alg.~\ref{alg:wcp} the error was estimated by the norm of vector $b_k = A^\t x_k.$ In the tensor algorithm we should estimate the Frobenius norm of the matrix~$\A \x_1 x_k$ instead, but it can not be evaluated by a small number of tenvecs and we substitute it by the spectral norm $\err \eqdef \|\A \x_1 x_k\|_2.$ 
To estimate this, we use the power method for $n_2 \times n_3$ matrix $\A \x_1 x_k^\t.$ It is initialized by some randomly chosen unit vector $z$ of size $n_3.$ 
In the power method, $\sigma$ converges to the maximum singular value of matrix $B_k = \A \x_1 x_k^\t$ and $\tilde B_k = \sigma y z^\t$ converges to the best rank-one approximation of $B_k.$
Since the high precision is not necessary for the error estimation, we can satisfy with the fixed small number $\ppow$ of power iteration steps.  
If only the fixed-rank stopping criteria is desired, power iterations on Steps 9-13 of Alg.~\ref{alg:wt1} can be omitted. 

As well as in the matrix case, we can meet with the \emph{breakdown} if new vector $x$ have almost zero (neglectable in the machine precision) component $x'$ orthogonal to the subspace~$X_{k-1}.$ We can try to fix it by repeating the current step with another selection of $y_k, z_k,$ or choose to terminate the algorithm. 

Mode-$2$ and mode-$3$ bases $V$ and $W$ can be computed by the same algorithm after obvious permutation of modes. 
Directly from the Statement~\ref{stat5} in matrix case we derive the following theorem, that proves the convergence of tensor methods based on Alg.~\ref{alg:wt1} for the exact-rank case.
\begin{theorem}\label{thm}
For the tensor $\A$ with mode sizes $n_1, n_2, n_3$ and ranks $r_1, r_2, r_3,$ three applications of Alg.~\ref{alg:wt1} return bases $U,$ $V$ and $W$ of sizes $n_1\x r_1,$ $n_2 \x r_2$ and $n_3 \x r_3$ that allow the Tucker decomposition $\A = \G \x_1 U \x_2 V \x_3 W$  with core $\G$ given by~\eqref{eq:core}, providing that computations are not terminated by breakdowns.
\end{theorem}
\begin{remark}
The convergence in the exact-rank case is guaranteed for any choice of the leading vectors $y_k$ and $z_k$ that does not lead to the breakdown. \end{remark}

In the section~\ref{LEAD} we propose a number of strategies for the selection of the leading vectors, that are based on the different optimization ideas and lead to the methods with different complexity and convergence properties. However, the convergence in the exact-rank case, which we consider as the necessary requirement, persists for all methods based on the Alg.~\ref{alg:wt1}.
\begin{remark}
The minimal Krylov recursion Alg.~\ref{alg:min} can be also considered as variant of Alg.~\ref{alg:wt1} with the special choice of leading vectors. 
\end{remark}
 
To achieve the reduction of mode-1 rank after the first iteration of the Alg.~\ref{alg:wt1} we should approximate the tensor by 
\begin{equation}\label{eq:rank1}
\A \approx \tilde \A_1 = x_1 \o B_1 \qquad \mbox{with} \qquad B_1 = \A \x_1 x_1.
\end{equation}
This does not comply the restrictions imposed to reach the linear complexity. 
Therefore, we consider Alg.~\ref{alg:wt1} as a method to generate an approximation of the dominant mode subspaces that is guaranteed to converge in $r$ steps for the tensor with the mode rank $r,$ but that does not reduce the mode ranks by one on each iteration.
\begin{remark}
The mode ranks of the rank-$(r_1,r_2,r_3)$ tensor in general can not be reduced by the elimination of the rank-$(1,1,1)$ approximation.
\end{remark}

To take use of the approximation~\eqref{eq:rank1}, we can further approximate $B_1$ by low-rank (or even rank-one) format. To do this by tenvec operations, we can use `internal' Wedderburn elimination steps for $B_1.$ 
This approach directly follows the basic idea of~\cite{ost-tucker-2008}, where the cross approximation method is generalized to $3$-tensors.    
Vectors  $x_k$  can be associated with the mode fibres selected in the tensor, and $B_k$ with the slices, for which the internal cross approximation scheme is used. As well as in the matrix case, we emphasize the very important link between the Krylov subspaces approach and cross approximation methods, that is established by the Wedderburn framework. 

In the following we restrict the discussion to only rank-one approximation of $B_k$ and pay by losing the mode rank reduction property, but result in the methods that are more efficient since they have `symmetric' behaviour in respect to different modes.

\subsection{Selecting leading vector}\label{LEAD}
By the Theorem~\ref{thm}, in Alg.~\ref{alg:wt1} the every choice of the leading vectors $y_k,  z_k$ that does not lead to the breakdown, ensures the convergence in the exact-rank case. 
However  the different choices of leading vectors result in the approximations with different accuracy until the exact representation is found.
Therefore, we should choose the leading vectors in a clever way, that ensures the fast convergence to the dominant subspaces and is computationally feasible. In the following we propose four strategies that lead to different maximization problems and result in different complexity estimates and convergence properties.

\subsubsection{SVD-like strategy (Wsvd)}
We can apply Alg.~\ref{alg:wt1} three times and compute $U, V, W$ in the completely independent processes, even using three processors on a distributed memory system.
In algorithm for $U = X_k$ the best way to keep from breakdowns is to choose $y_k \o z_k$ that maximizes the orthogonal component $\|x'\|.$ 
\begin{equation}\label{eq:wsvd}
\begin{split}
      y_k, z_k & = \arg \max_{\|y\|=\|z\|=1}\left\| (I - X_{k-1} X_{k-1}^\t) (\A yz) \right\|
                   = \arg \max_{\|y\|=\|z\|=1}\left\| \B  yz \right\|, \\
      \B &= \A \x_1 (I - X_{k-1} X_{k-1}^\t).
\end{split} 
\end{equation}
This is the direct analogy with the SVD approach~\eqref{eq:best} in the matrix case. 
The difference is that resulted $x_k$ and $y_k \o z_k$ are not singular vectors of the unfolding $A = A^{(1)},$ since the maximization is done w.r.t. `long' right vectors with tensor product structure. This is the reason why the best tensor rank-$(r,r,r)$ approximation can not be computed by the simultaneous elimination of the best rank-$(1,1,1)$ approximations.

Nevertheless with this choice we are safe from breakdowns. 
Zero-valued orthogonal component $\|x'\|$ appears only if $\max_{\|y\|=\|z\|=1}\left\| \B  yz \right\| = 0$ and hence $\B=0,$ 
which means that the \emph{exact} representation $\tilde \A_k \eqdef \A \x_1 (X_kX_k^\t)$ with mode-$1$ rank $r_1=k$ is computed for $\A.$ 
The following theorem shows that `machine precision breakdown' also happens only when the approximation is of the machine precision accuracy.
\begin{theorem}
If the breakdown $\|x'\| < \tol \|x\|$ is met on the step $k+1$ of the Alg.~\ref{alg:wt1}, than the computed subspace $X_k$ provides  the approximation $\tilde \A_k = \A \x_1 (X_k X_k^\t)$  with the Tucker rank $r_1 = k$ and accuracy $\|\A - \tilde \A_k\|_\tn < \tol\|\A\|_\tn$ in spectral norm. 
\end{theorem}
\begin{proof}
On the step $k+1$ we choose $y_{k+1},z_{k+1} = \arg\max_{\|y\|=\|z\|=1} \|\B yz \|.$ 
The norm of the orthogonal component reads
\begin{equation}\nonumber
\begin{split}
\|x'\| =  & {} \|\B y_{k+1} z_{k+1}\| =  \max_{\|y\|=\|z\|=1} \|\B  yz\| {} \\ {} = &  {}
           \max_{\|y\|=\|z\|=1}\max_{\|x\|=1} (\B  yz, x) =  \max_{\|x\|=\|y\|=\|z\|=1} \f{\B}{x\o y\o z} \eqdef \|\B\|_\tn.
\end{split}
\end{equation}
Also for $x = \A y_{k+1} z_{k+1}$ it holds
$$
\|x\| = \|\A y_{k+1} z_{k+1}\| = \max_{\|x\|=1} (\A y_{k+1} z_{k+1}, x) \leq \max_{\|x\|=\|y\|=\|z\|=1} \f{\A}{x \o y \o z} \eqdef \|\A\|_\tn.
$$
Now the breakdown criteria $\|x'\| < \tol \|x\|$ gives $\|\B\|_\tn <\tol \|\A\|_\tn.$ Finally we  write 
$$
\B = \A \x_1 (I - X_k X_k^\t) = \A - \A \x_1 (X_k X_k^\t) = \A - \tilde \A_k,
$$
that completes the proof since mode-$1$ rank of $\tilde A_k$ is equal to $r_1  = \rank X_k = k.$
\end{proof}

\begin{corollary}
Alg.~\ref{alg:wt1} with the SVD-like strategy~\eqref{eq:wsvd} applied to tensor $\A$ with mode ranks $r_1, r_2, r_3$ computes bases $U, V, W$ that allow the exact representation $\A = \G \x_1 U \x_2 V \x_3 W$ after $r_1, r_2$ and $r_3$ iterations, respectively.
\end{corollary}

Note that since $X_k = [X_{k-1} \: x_k]$ is orthogonal, it holds
$$
\A \x_1 x_k^\t = \A \x_1 \left( (I-X_{k-1} X_{k-1}^\t) x_k \right)^\t = \left( \A \x_1 (I - X_{k-1}X_{k-1}^\t) \right) \x_1 x_k^\t = 
\B \x_1 x_k^\t.
$$
The maximization problem~\eqref{eq:wsvd} can be solved by the ALS Alg.~\ref{alg:als} applied for $\B.$
In this case power iterations for matrix $\A \x_1 x_k^\t$ on Step 9-13 of Alg.~\ref{alg:wt1} can be omitted, since the error estimate $\err = \|\A \x_1 x_k^\t\|_2 = \|\B \x_1 x_k^\t\|_2 $ is actually computed in the ALS iterations.

On the step $k$ of the Wedderburn method inner ALS iterations cost $3 \pals$ tenvecs and $\O(\pals nk)$ operations for the orthogonalization. 
This summarizes to $3 \pals r$ tenvecs and $\O(\pals nr^2)$ additional operations for one dominant subspace.

\subsubsection{Lanczos-like strategy (Wlnc)} 
We can also use the analogy with the Lanczos choice~\eqref{eq:choice} by taking unit $y_k \o z_k \approx \A \x_1 x_k^\t / \|\A \x_1 x_k^\t\|.$ 
This leads to the dual maximization problem
\begin{equation}\label{eq:wlnc}
\begin{split}
      y_k, z_k &= \arg \max_{\|y\|=\|z\|=1}\left\|  \A \x_1 x_k^\t \x_2 y^\t \x_3 z^\t \right\|
            = \arg \max_{\|y\|=\|z\|=1}\left\| y^\t B z \right\| \\
      B &   = \A \x_1 x_k^\t.
\end{split} 
\end{equation}
The solution can be accomplished by $\ppow$ steps of the power iteration method applied to matrix $B.$ 
On the step $k$ of the Wedderburn method it costs $2 \ppow$ tenvecs. 
Note that the maximization problem~\eqref{eq:wlnc} is actually solved on Steps 9-13 of Alg.~\ref{alg:wt1} by power iterations that computes best rank-$1$ approximation of $\A\x_1 x_k^\t$ to estimate the norm of the residual. 
Thus, the Wlnc pivoting strategy requires only to set $y_k := y$ and $z_k:=z$ after power iterations as new vectors of the Wedderburn  process.

\subsubsection{Restricted SVD-like strategy (WsvdR)} 
Three Wedderburn elimination algorithms can be used to extend all mode bases  \emph{simultaneously}.
Suppose $k-1, l$ and $m$ steps were done to compute mode subspaces $X_{k-1}, Y_l, Z_m.$
Then on the step $k$ of the Wedderburn process for the mode-$1$ subspace we can make use of $V = Y_l$ and $W = Z_m$ by \emph{restricting} the maximization~\eqref{eq:wsvd} to the tensor product of these subspaces. 
Therefore, we take $y_k=Y_l\yy_k, z_k = Z_m\zz_k$  and solve
\begin{equation}\label{eq:wsvdr}
\begin{split}
      \yy_k, \zz_k & = \arg \max_{\|\yy\|=\|\zz\|=1}\left\| (I - X_{k-1} X_{k-1}^\t) \left( \A (Y_l\yy)(Z_m\zz) \right) \right\|
                   = \arg \max_{\|\yy\|=\|\zz\|=1}\left\| \B  \yy \zz \right\| \\
      \B &= \A \x_1 (I - X_{k-1} X_{k-1}^\t) \x_2 Y_l^\t \x_3 Z_m^\t.
\end{split} 
\end{equation}
This  maximization problems exactly matches the~\eqref{eq:optmax}!
\begin{remark}
The optimized MKR is a variant of the Wedderburn process for tensors with restricted SVD-like strategy of pivoting.
\end{remark}

The WsvdR approach can lead to the slow convergence or breakdowns at the first iterations, but when $X_k, Y_l, Z_m$ become larger, chances to meet the breakdown vanish, since~\eqref{eq:wsvdr} becomes close to the unrestricted maximization~\eqref{eq:wsvd}. 
Numerical experiments provided in the section~\ref{NUMc} show that sometimes the restricted SVD strategy gives even better results than unrestricted SVD strategy. The reason is probably that with this restrictions the  ALS iterations are not caught in the local minima.

Complexity on the step $k$ for the mode-$1$ subspace is $3\pals$ tenvecs and $\O(\pals nk)$ operations for the orthogonalization, total complexity is $9 \pals r$ tenvecs and $\O(\pals nr^2)$ additional operations.

\subsubsection{Restricted Lanczos-like strategy (WlncR)}
\begin{algorithm}[t]
\caption{Wedderburn with restricted Lanczos-like pivoting strategy (WlncR)}\label{alg:wt3}
\begin{algorithmic}[1]
\REQUIRE Tenvec subroutine for the tensor $\A,$  tolerance  $\tol,$ accuracy $\eps$
\ENSURE Approximation $\tilde A = \G \x_1 U \x_2 V \x_3 W$ such that $\|\A-\tilde\A\|_F \lesssim \|\tilde\A\|_F$
\item[\textbf{Initialization:}] Unit vectors $u, v, w,$ $\updx=\updy=\updz=\TRUE.  $ 
\STATE $x_1=\A vw/\|\A vw\|; \quad y_1=\A wu/\|\A wu\|; \quad z_1=\A uv/\|\A uv\|, \quad k=l=m=1$
\STATE  $X_1 = [x_1], Y_1=[y_1], Z_1=[z_1], \quad \G = \A \x_1 x_1^\t \x_2 y_1^\t \x_3 z_1^\t; \quad \nrm=\|\G\|_F$
\WHILE{$\updx$ or $\updy$ or  $\updz$}
      \item[{\small{x:}}] \COMMENT{Proceed with new vector $x$ if required}
      \IF{$\updu$}
      \STATE            For $B=\G(k,:,:)$ solve $B \approx: \tilde B = b \yy_k \zz_k^\t$ \COMMENT{best rank-one approximation}
      \STATE            $y_k := Y_l\yy_k, \quad z_k := Z_m\zz_k; \quad x := \A y_k z_k; \quad x' := (I-X_k X_k^\t)x$
                  \IF[Breakdown]{$\|x'\| < \tol \|x\|$}
                      \STATE   $\updx:=\FALSE$ \COMMENT{or repeat current iteration step with another $y_k, z_k$} 
                  \ELSE
                      \STATE   $x_{k+1}:=x'/\|x'\|; \quad X_{k+1} := [X_k \: x_{k+1}] $ 
                      \STATE   Enlarge $\G$ by $\G(k+1,:,:) := \A \x_1 x_{k+1}^\t \x_2 Y_l^\t \x_3 Z_m^\t$
                      \STATE  $\err:=\|\G(k+1,:,:)\|_F, \quad \nrm^2 := \nrm^2 + \err^2, \quad k:=k+1$
                      \IF[Convergence]{$\err < \eps \nrm$}
                              \STATE $\updx:=\FALSE$
                      \ENDIF
                  \ENDIF
      \ENDIF
      \item[{\small{y:}}] \COMMENT{Proceed with new vector $y$ if required}
      \item[{\small{z:}}] \COMMENT{Proceed with new vector $z$ if required}
\ENDWHILE
\RETURN $U = X_k, V = Y_l, W = Z_m, \quad \tilde A = \G \x_1 U \x_2 V \x_3 W$
\end{algorithmic}
\end{algorithm}

Finally, we combine the Lanczos-like selection of leading vectors~\eqref{eq:wlnc} and the restricted maximization.
\begin{equation}\label{eq:wlncr}
\begin{split}
      \yy_k, \zz_k & = \arg \max_{\|\yy\|=\|\zz\|=1}\left\|  (Y_l\yy)^\t (\A\x_1 x_k^\t) (Z_m\zz) \right\|
                   = \arg \max_{\|\yy\|=\|\zz\|=1}\left\| \yy^\t B \zz  \right\| \\
      B &= \A \x_1 x_k^\t \x_2 Y_l^\t \x_3 Z_m^\t, \qquad y_k=Y_l\yy_k, \quad z_k = Z_m\zz_k.
\end{split} 
\end{equation}
If only bases  $U, V, W$ for the approximation are required, the maximization can be accomplished by $\ppow$ steps of power iteration method applied to the matrix $B.$ At every step of the Wedderburn method it requires $2 \ppow$ tenvecs.
But if the core tensor is also desired, the efficiency can be highly improved by precomputing $B$ as $l \x m$ matrix.
Comparing~\eqref{eq:wlncr} and~\eqref{eq:core} we note that $B$ is exactly the last mode-$1$ slice from the optimal core for the approximation of $\A$ in bases $U = X_k, V = Y_l, W=Z_m$
$$
\G(k,:,:) = \A \x_1 (U(:,k))^\t \x_2 V^\t \x_3 W^\t = \A \x_1 x_k^\t \x_2 Y_l^\t \x_3 Z_m^\t = B.
$$
Thus if we need the core $\G,$ we prefer to compute it slice-by-slice in the Wedderburn process. 
Then we can apply standard matrix tools (SVD or cross methods) to find the best rank-one approximation of $B$ and solve~\eqref{eq:wlncr} without additional tensor operations.
This version of the Wedderburn elimination algorithm for tensors is summarized in the Alg.~\ref{alg:wt3}.

\subsection{Comparison of the algorithms}
\begin{table}[tb]
\caption{Complexity of the algorithms for rank-$(r_1, r_2, r_3)$ tensor approximation} \label{tab1}
\begin{center}
\begin{tabular}{ccc|c}
name  & description                                                            & output         & tenvecs         \\ \hline
MKR   & Min. Krylov recursion~\cite{eldensavas-krylov-2009} Alg.~\ref{alg:min} & $U, V, W$      & $3r$            \\
Wsvd  & Alg.~\ref{alg:wt1} with strategy~\eqref{eq:wsvd}                       & $U, V, W$      & $9\pals r + 3r$ \\ 
Wlnc  & Alg.~\ref{alg:wt1} with strategy~\eqref{eq:wlnc}                       & $U, V, W$      & $6\ppow r + 3r$ \\ 
WsvdR & Alg.~\ref{alg:wt1} with strategy~\eqref{eq:wsvdr}, Alg.~\ref{alg:opt}  & $U, V, W$      & $9\pals r + 3r$ \\
WlncR & Alg.~\ref{alg:wt1} with strategy~\eqref{eq:wlncr}, Alg.~\ref{alg:wt3}  & $U, V, W$      & $6\ppow r + 3r$ \\
WlncR &       --- // ---                                                       & $U, V, W, \G$  & $r^2 + 3r$
\end{tabular}
\end{center}
\end{table}
To compare the versions of the proposed algorithm, we give their complexities in Table~\ref{tab1}. 
If only subspaces $U, V, W$ for Tucker rank-$(r_1,r_2,r_3)$ approximation are required, than the minimal Krylov recursion~\cite{eldensavas-krylov-2009} is fastest in terms of number of tenvecs used, since it requires only $3r$ tensor operations. 
All versions of the Wedderburn elimination Alg.~\ref{alg:wt1} also require $\O(r)$  tenvecs, with factor depending only on the selected number of ALS or power iterations. 
Note that for $\pals = \ppow,$  Lanczos-like pivoting strategy and SVD-like pivoting take roughly the same time, but SVD-like pivoting is guaranteed to be free from breakdowns. This differ tensor algorithms from the matrix case.

If the core $\G$ for Tucker approximation is required, we generally need additional $r^2$ tenvecs to evaluate it by~\eqref{eq:core}. 
In this case, the fastest version of the Wedderburn elimination algorithm is the WlncR Alg.~\ref{alg:wt3}. 
It uses \emph{exactly} the same number of tenvecs as the MKR, but shows better convergence in numerical experiments.

\section{Numerical examples}\label{NUM}
The numerical experiments presented in this section were performed on the Intel Xeon Quad-Core E5504 CPU running at $2.00$~GHz. 
In the section \ref{NUMs} we use MATLAB version 7.7.0, in sections~\ref{NUMc},~\ref{NUMt2t} we use Intel Fortran compiler version 11.1 and BLAS/LAPACK routines provided by the MKL library.

\subsection{Sparse tensors}\label{NUMs}
The tensor decomposition of a sparse large-dimensional arrays is an important tool in the network analysis, that is widely used now in the information science, sociology and many other disciplines. 
We consider one example from~\cite{facebook}, where a very nice introduction to network analysis is given. 

\begin{figure}[t]
\caption{Approximation accuracy for the Caltech Facebook graph}
\begin{center}
\includegraphics[width=.65\textwidth]{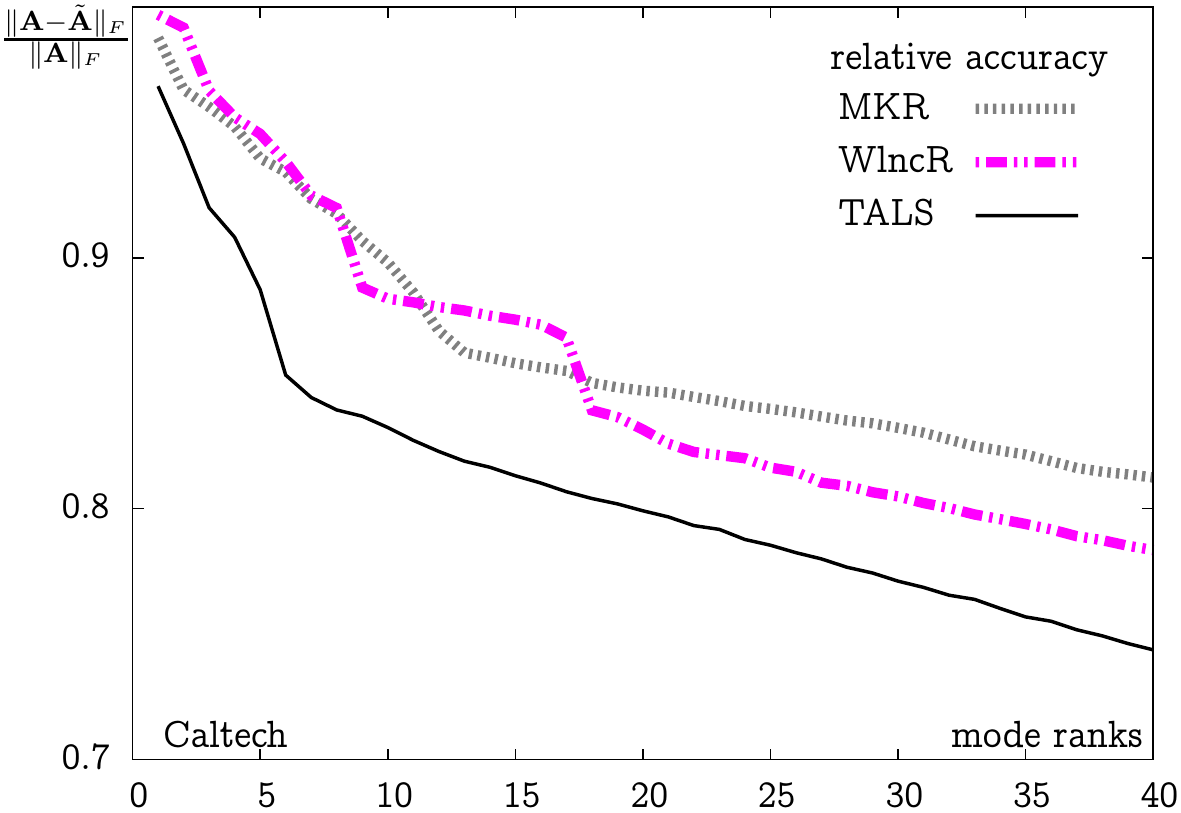} 
\end{center}
\cntr{
We approximate the $597 \x 597 \x 8^2$ tensor using three algorithms and
plot relative accuracy in the Frobenius norm for different values of mode ranks.
TALS is the Tucker-ALS, MKR is the minimal Krylov recursion, WlncR is the Wedderburn elimination with the Lanczos-like restricted pivoting. 
} 
\label{fig:caltech}
\end{figure}

The example is the graph of the Facebook social network from Caltech University, coupled with the dormitory information. 
The graph represents relations between $n=597$ people living in $m=8$ dorms by sparse $n \x n \x m \x m$ $4$-tensor $\A.$ 
Every entry $a_{ijpq}=1$ means that the person number $i$ lives in the dorm number $p$ and links to the person number $j,$ who lives in the dorm $q.$ There are $25646$ unit elements in $\A=[a_{ijpq}],$ and all others are zeroes.
Since the discussed methods apply to $3$-tensors, we join dorm indices $p$ and $q$ in multi-index and consider $\A=[a_{i,j,pq}]$ as the $3$-tensor of size $n \x n \x m^2.$ 
The relative accuracy of the approximation in the Frobenius norm is shown on Fig.~\ref{fig:caltech}. 
We compare three methods: the Tucker-ALS \cite{tuckerals-1980,lathauwer-rank1-2000} (MATLAB implementation from Tensor Toolbox~\cite{kolda-toolbox-2007}), the minimal Krylov recursion (Alg.~\ref{alg:min}) and the Wedderburn elimination method with Lanczos-like restricted pivoting (Alg.~\ref{alg:wt3}). 
We note that the WlncR Alg.~\ref{alg:min} converges slowly on the first iterations, when dimensions of accumulated bases are small, imposing serious restrictions in the maximization~\eqref{eq:wlncr}. 
However, for the larger values of ranks, the WlncR becomes more accurate than the MKR method. 

The Tucker-ALS is the most accurate but computationally demanding method. 
On each Tucker-ALS iteration the Tucker factors are subsequently updated as follows. 
With two factors fixed, for example $V$ and $W,$ we compute the  $n_1 \x r_2 \x r_3$ tensor
\begin{equation}\label{eq:tals}
\B = \A \x_2 V^\t \x_3 W^\t, \qquad B = B^{(1)}
\end{equation}
Then for the $n_1 \x r_2r_3$ unfolding $B = B^{(1)}$ we find the approximation  $B \approx \tilde B =: U G$ with the $n_1 \x r_1$ orthogonal matrix $U$ that is a new Tucker factor and the $r_1 \x r_2 r_3$ matrix $G$ that is then reshaped to new $r_1 \x r_2 \x r_3$ core tensor $\G.$ 
The rank-$r_1$ approximation of $B$ can be computed by the SVD or cross approximation methods (c.f.~\cite{tee-cross-2000}). 
The evaluation of~\eqref{eq:tals} requires $r_2 r_3$ tenvecs, which results in $3 r^2$ tenvecs on each iteration, much more than the complexity of the MKR and the WlncR.

We do not provide timings in this section, because MATLAB-based computations are usually far from being highly optimized. 
Some timings will appear in the next sections for the Fortran implementation of discussed algorithms.

This example was introduced to us by Prof. Lars Eld\'en. In~\cite{eldensavas-krylov-2010} more experiments with the approximation of sparse tensors by different algorithms are provided and results are compared for the truncated HOSVD, the minimal Krylov recursion, several modifications of the MKR and the WlncR Alg.~\ref{alg:wt3} that was implemented by authors of~\cite{eldensavas-krylov-2010} with minor modifications.

\subsection{Compression from canonical to Tucker format}\label{NUMc}

Multidimensional data often appear in modern modelling programs in canonical form $(C).$ 
For example, in chemical modelling programs, e.g. PC GAMESS and MOLPRO, the \emph{electron density function} is given as a sum of the tensor product of one-dimensional Gaussians.
However, even for simple molecules, the number of terms in the decomposition obtained by MOLPRO may be too large for practically feasible computations. 
In order to make computations efficient, the further approximation (recompression) to the Tucker format can be performed.  
The accuracy of the desired approximation can vary in different applications. 
For some quantum chemistry problems the very precise approximation (with ten or more significant digits) is required.

Such recompression was done in~\cite{mpi-chem3d-2009} using the Tucker-ALS algorithm~\cite{tuckerals-1980,lathauwer-rank1-2000}, in \cite{khor-ml-2009} by the Tucker-ALS with the initial guess obtained from the coarser grids, in~\cite{fkst-chem-2008} by the Cross3D algorithm~\cite{ost-tucker-2008}, in~\cite{ost-chem-2009} by the individual cross approximation of canonical factors and in~\cite{sav-rr-2009} by the cross approximation of Gram matrices of unfoldings. 
For an $\nnn$ tensor given in canonical form with $R$ terms, each tenvec costs $3nR$ operations, and the proposed versions of Alg.~\ref{alg:wt1} can be applied to compute the Tucker approximation efficiently, even for large $n$ and $R.$ 

\begin{figure}[p]
\caption{Approximation accuracy of the methane electron density}
\begin{center}
\includegraphics[width=.45\textwidth]{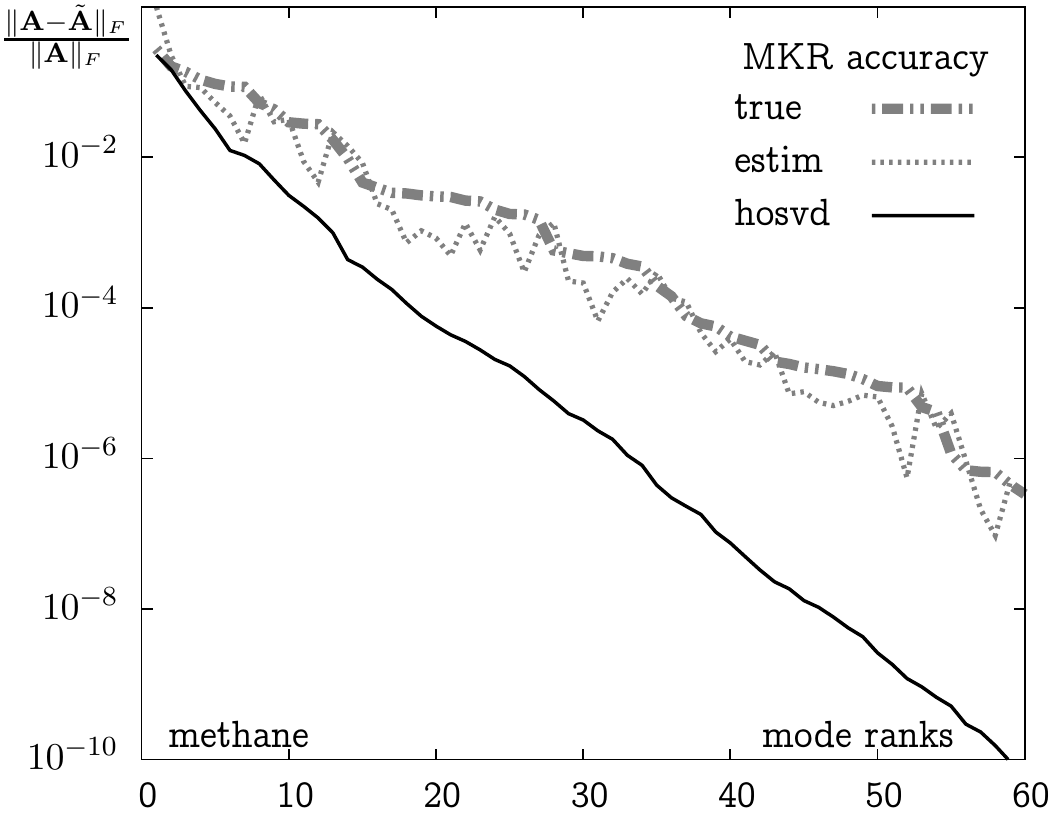} 
\end{center}
\begin{center}\hfil
\includegraphics[width=.45\textwidth]{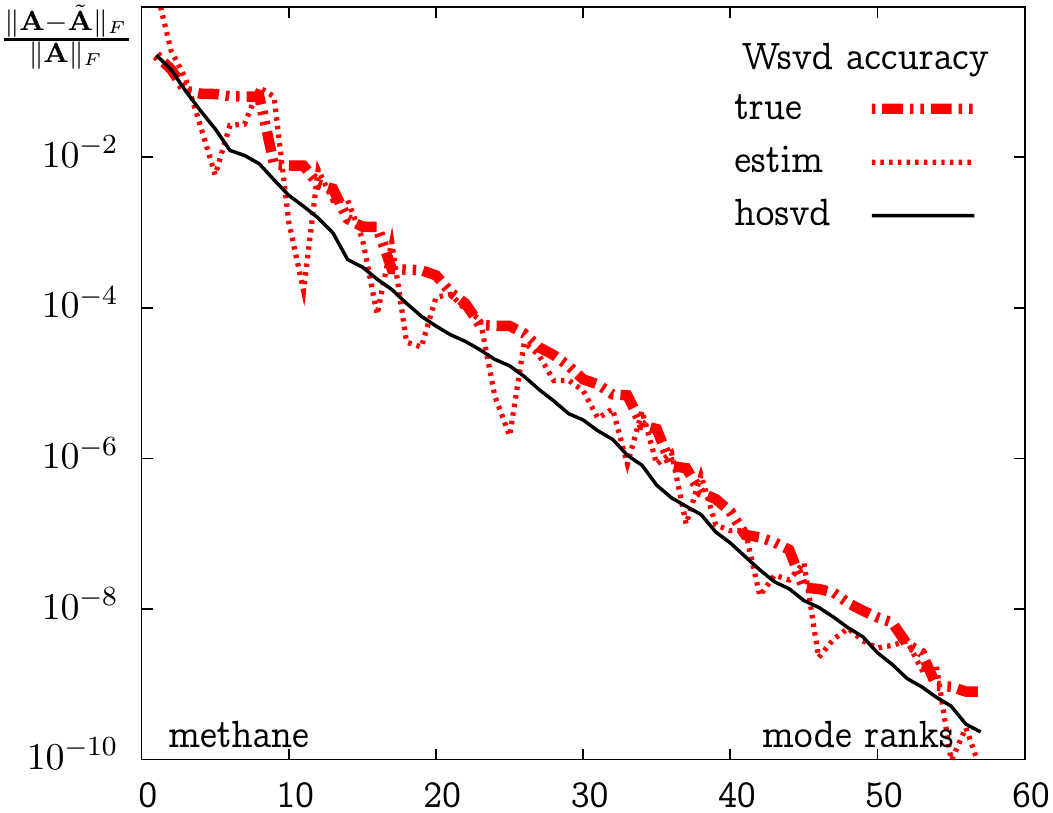} \hfil
\includegraphics[width=.45\textwidth]{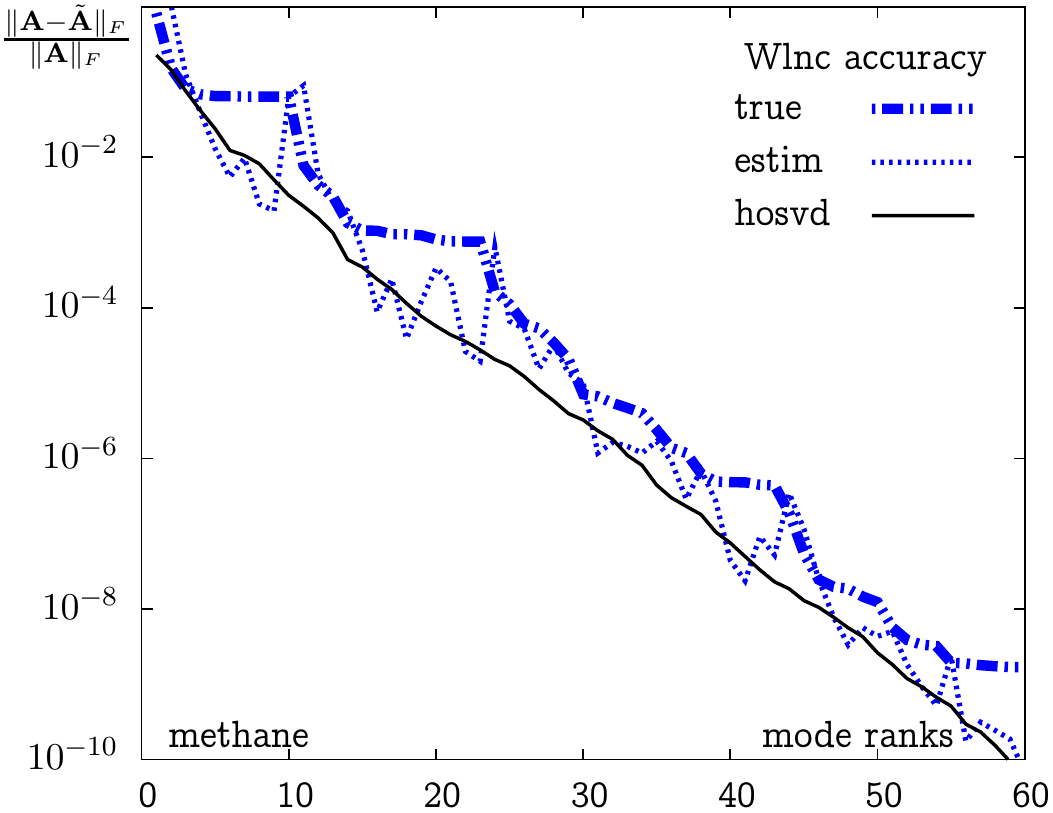} \hfil
\end{center}
\begin{center}\hfil
\includegraphics[width=.45\textwidth]{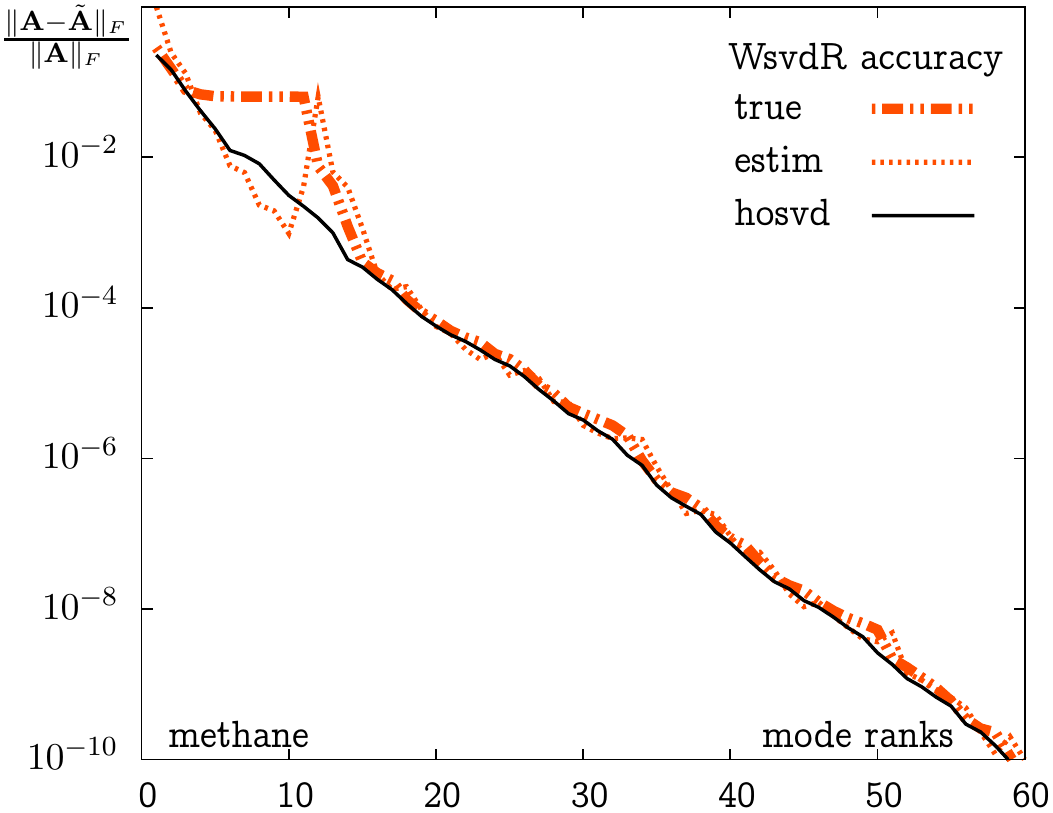} \hfil
\includegraphics[width=.45\textwidth]{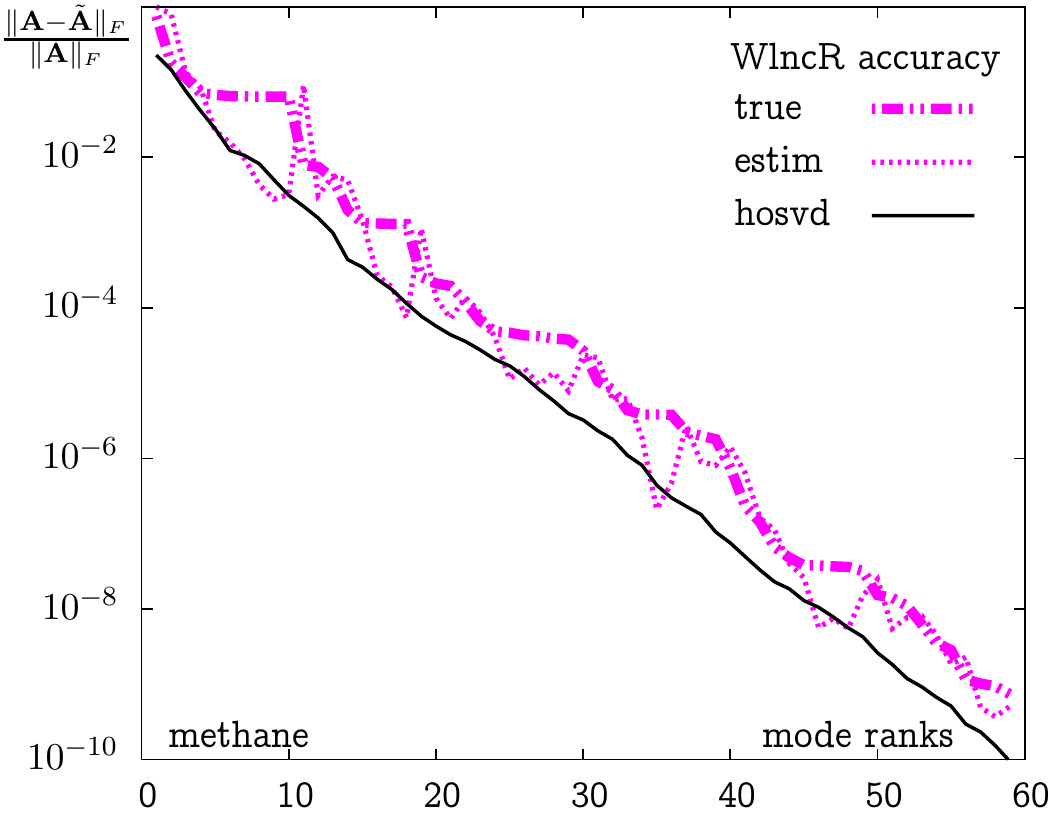} \hfil
\end{center}
\cntr{For $5121 \x 5121 \x 5121$ tensor given in the canonical form with rank $1334,$ we compute the Tucker approximation by algorithms, listed in Tab.~\ref{tab1}, and plot the relative accuracy in the Frobenius norm for different values of mode ranks. 
On each graph the thick dashed line shows the true error of approximation, the thin dashed line shows the estimate of error and the thin solid line shows the accuracy of the HOSVD approximation for the reference. 
} \label{fig:c1}
\end{figure}

\begin{figure}[p]
\caption{Approximation accuracy of the glycine electron density}
\begin{center}
\includegraphics[width=.45\textwidth]{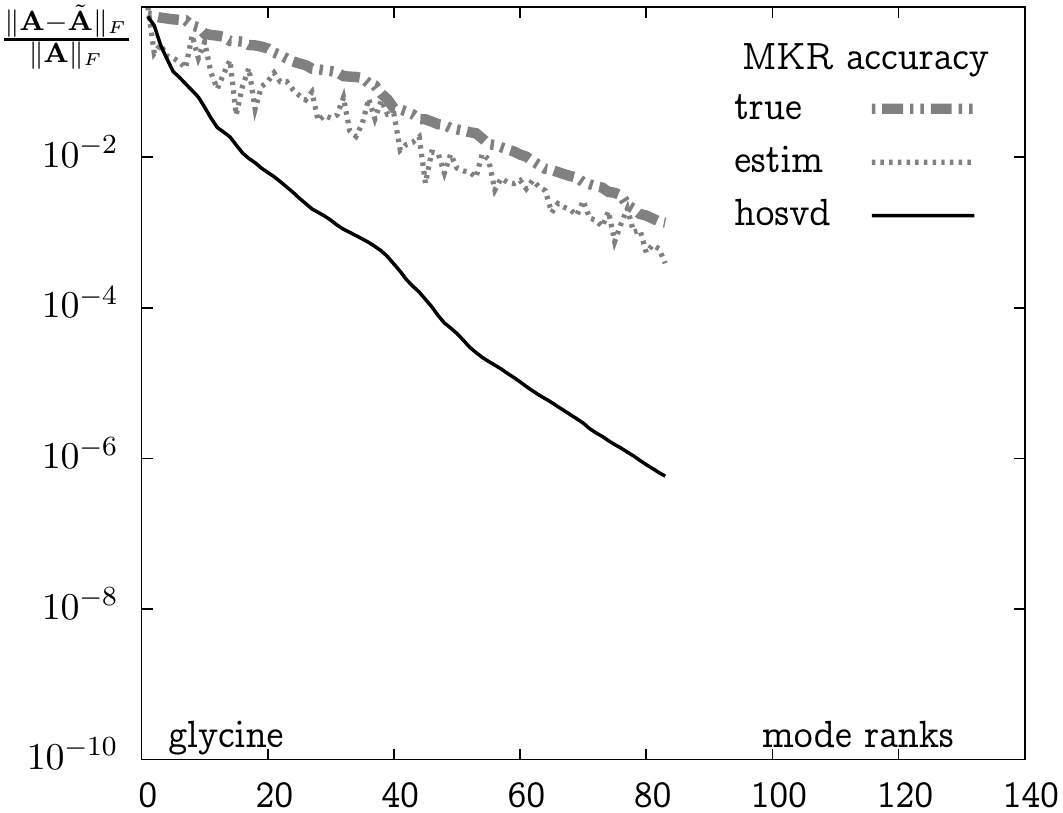} 
\end{center}
\begin{center}\hfil
\includegraphics[width=.45\textwidth]{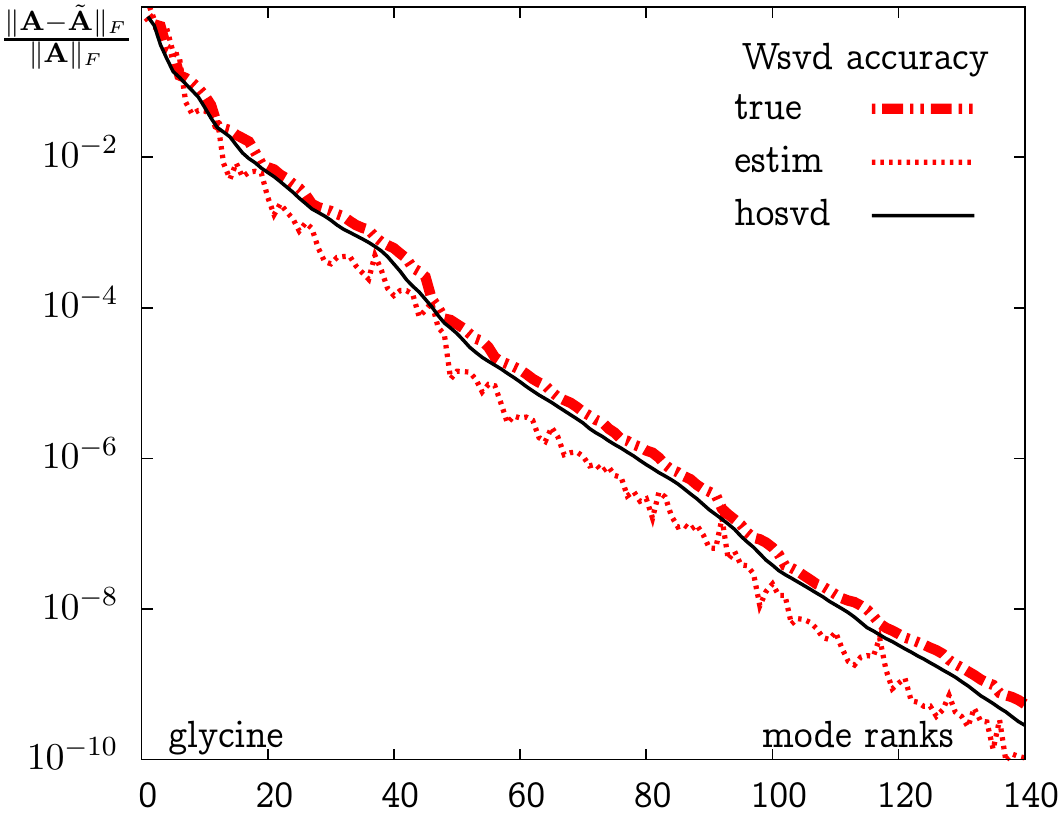} \hfil
\includegraphics[width=.45\textwidth]{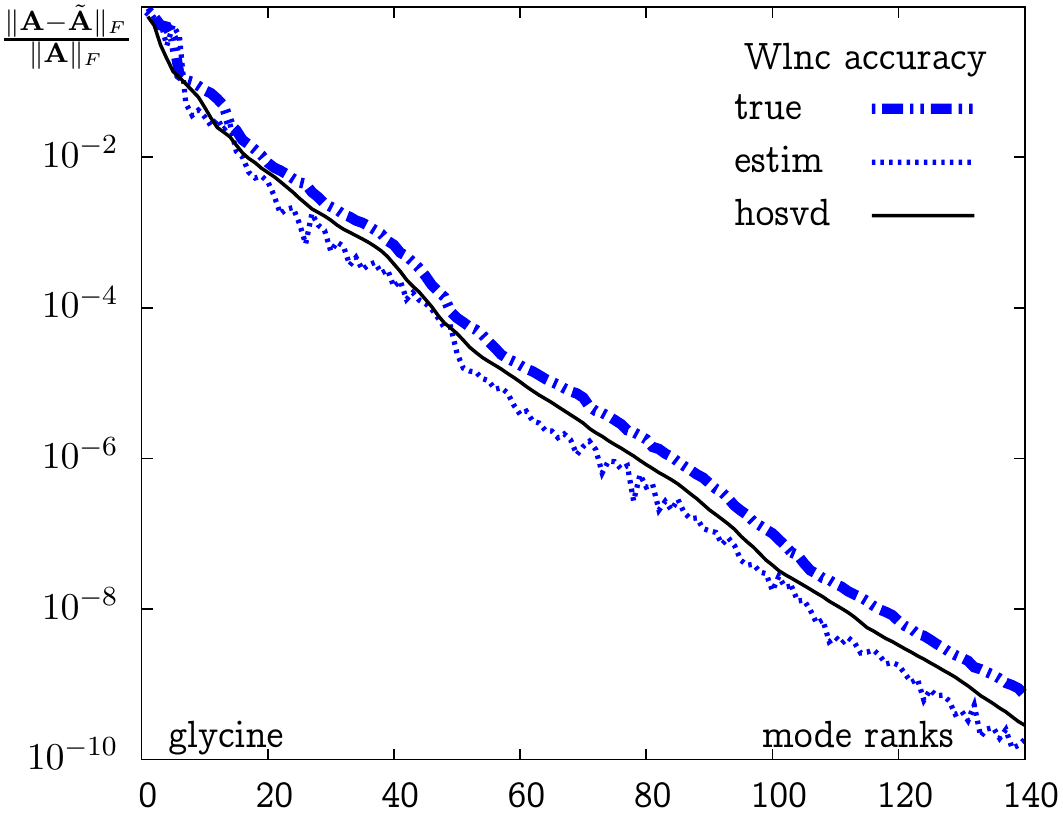} \hfil
\end{center}
\begin{center}\hfil
\includegraphics[width=.45\textwidth]{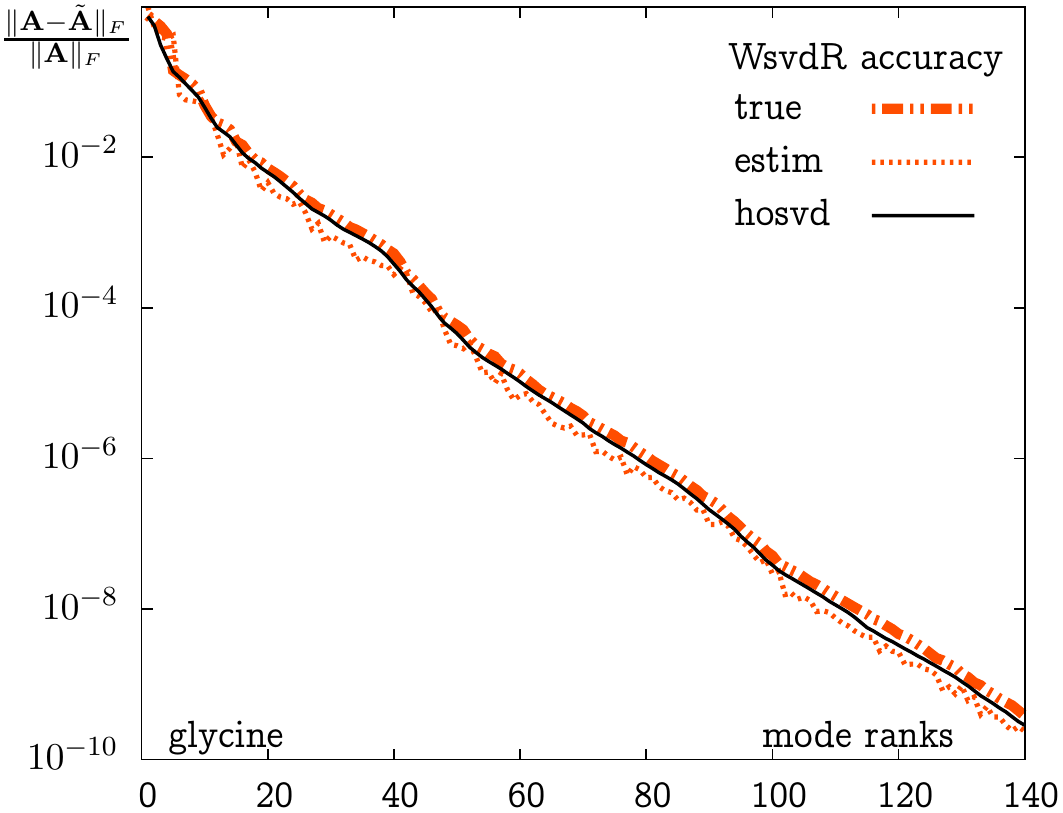} \hfil
\includegraphics[width=.45\textwidth]{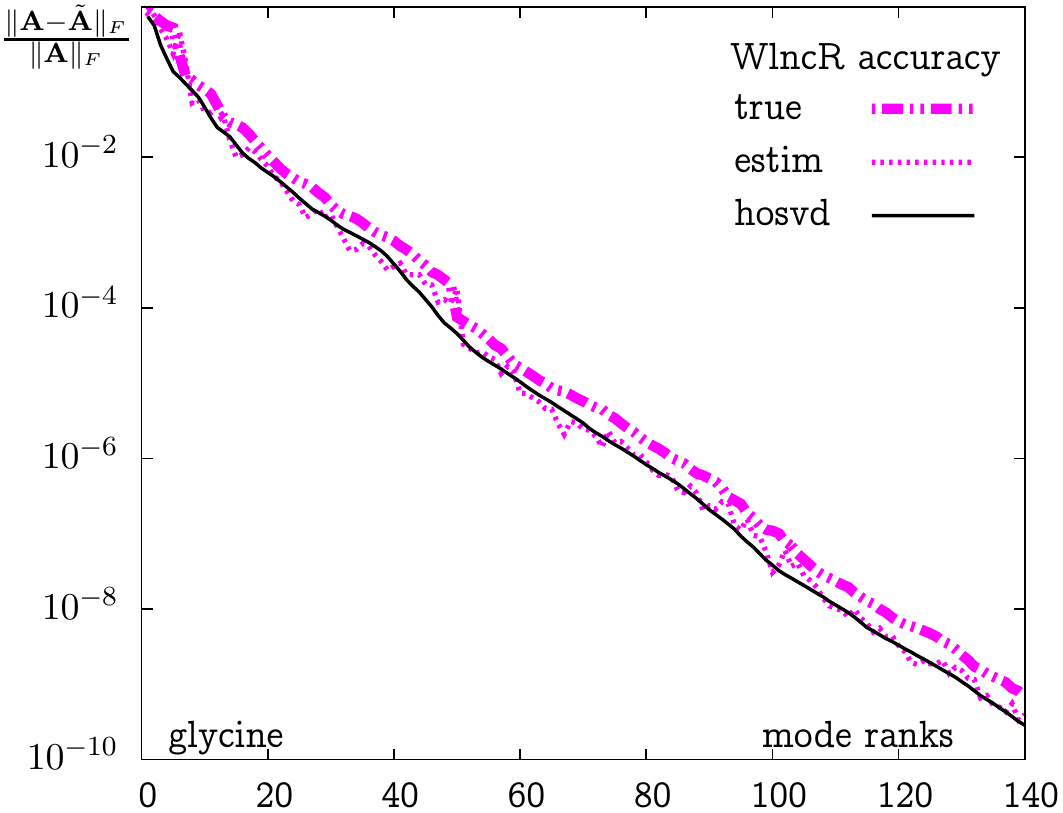} \hfil
\end{center}
\cntr{For $5121 \x 5121 \x 5121$ tensor given in the canonical form with rank $9208,$ we compute the Tucker approximation by algorithms, listed in Tab.~\ref{tab1}, and plot the relative accuracy in the Frobenius norm for different values of mode ranks. 
On each graph the thick dashed line shows the true error of approximation, the thin dashed line shows the estimate of error and the thin solid line shows the accuracy of the HOSVD approximation for the reference. 
} \label{fig:c2}
\end{figure}

\begin{table}[t]
\caption{Time for approximation of electron density}
\begin{center}
\begin{tabular}{cc|c|cccc|cc}
molecule    & accuracy   & MKR   & Wsvd   & Wlnc   & WsvdR  & WlncR & cross  & TALS($1$)\\ \hline
methane     & $10^{-4}$  & $1.0$ & $6.1$  & $4.4$  & $2.0$  & $0.5$ & $1.0$  & $2.6$    \\
$R=1334$    & $10^{-6}$  & $1.7$ & $10.7$ & $8.2$  & $3.7$  & $0.8$ & $1.4$  & $9.6$   \\
            & $10^{-8}$  &  ---  & $14.5$ & $11.1$ & $5.1$  & $1.4$ & $1.9$  & $30$   \\
            & $10^{-10}$ &  ---  & $20.4$ & $16.1$ & $6.9$  & $2.0$ & $2.9$  & $59$   \\ \hline
ethane      & $10^{-4}$  & $2.7$ & $17$   & $15$   & $6.7$  & $1.6$ & $3.0$  & $4.6$    \\
$R=3744$    & $10^{-6}$  & $5.2$ & $32$   & $25$   & $12$   & $2.8$ & $4.9$  & $17$    \\
            & $10^{-8}$  & ---   & $45$   & $35$   & $17$   & $3.9$ & $6.3$  & $42$    \\
            & $10^{-10}$ & ---   & $61$   & $50$   & $23$   & $5.3$ & $8.2$  & $83$    \\ \hline
ethanol     & $10^{-4}$  & $8.0$ & $54$   & $45$   & $17$   & $4.7$ & $8.1$  & $22$   \\
$R=6945$    & $10^{-6}$  & $14$  & $89$   & $74$   & $30$   & $8.4$ & $13$   & $81$   \\
            & $10^{-8}$  &  ---  & $135$  & $108$  & $45$   & $13$  & $17$   & $194$    \\
            & $10^{-10}$ &  ---  & $180$  & $145$  & $61$   & $18$  & $22$   & $391$    \\ \hline
glycine     & $10^{-4}$  &  ---  & $85$   & $69$   & $37$   & $7.5$ & $24$   & $32$   \\
$R=9208$    & $10^{-6}$  &  ---  & $131$  & $112$  & $57$   & $13$  & $33$   & $96$   \\
            & $10^{-8}$  &  ---  & $200$  & $160$  & $80$   & $18$  & $43$   & $211$    \\
            & $10^{-10}$ &  ---  & $268$  & $211$  & $114$  & $24$  & $60$   & $412$    \\
\end{tabular}
\end{center}
\cntr{For $5121 \x 5121 \x 5121$ tensor given in the canonical form with rank $R,$ we compute the Tucker approximation with different relative accuracy bound by different algorithms and show time in seconds. Algorithms MKR, Wsvd, Wlnc, WsvdR, WlncR are listed in Tab.~\ref{tab1}, 
`cross' algorithm is based on individual cross approximation of canonical factors~\cite{ost-chem-2009},
`TALS(1)' is time for one iteration of the Tucker-ALS method~\cite{tuckerals-1980,lathauwer-rank1-2000}
} \label{tab2}
\end{table}

We apply the discussed algorithm for the Tucker approximation of the electron density of some simple molecules, discretized on the uniform $\nnn$ tensor grid with $n=5121.$ 
The convergence of algorithms, i.e. the accuracy of rank-$(r,r,r)$ approximation for different $r$ is shown on Figs.~\ref{fig:c1},~\ref{fig:c2}. 
On each graph we compare the accuracy of the approximation computed by the Tucker-ALS \cite{tuckerals-1980,lathauwer-rank1-2000} (thin solid line) with the accuracy of the certain approximation method.
For each method, the internal estimate of the error $\err$ is shown by thin dashed line and real accuracy $\|\A - \tilde \A\|_F$ of the algorithm is shown by the thick dashed line.
The core of the Tucker approximation was computed by~\eqref{eq:core}. 
Since the evaluation of the whole $\nnn$ array for $n=5121$ requires one terabyte of memory and a lot of computational resources, we verify the accuracy of algorithms by comparison of the result with the Tucker approximation computed by the individual cross approximation of canonical factors~\cite{ost-chem-2009} with the accuracy set to $\eps=10^{-12}.$ 
The verification of the latter was done in~\cite{ost-chem-2009} by  the exhaustive verification on cluster platforms.
The residual between two tensors in the Tucker format is computed as proposed in~\cite{ost-latensor-2009}.

We note the slow convergence of the MKR for the methane electron density and the breakdown for glycine. 
All versions of the Wedderburn elimination converge much better, and for the larger glycine molecule the convergence is even more regular, than for methane. 
Accuracy of methods with the Lanczos-like pivoting is close to the optimal one, and accuracy of method with the restricted SVD-like pivoting is almost equal to optimal, except for first steps of the process, when  accumulated subspaces are small and impose significant restrictions on the selection of leading vectors.
Also, the internal error value, estimated by norm of $\A \x_1 x_k^\t$ (and the same for other modes), is less regular than the real accuracy, that decays monotonically. 
In the methods with unrestricted pivoting (see Wlnc) the internal error value is computed by the spectral norm of matrix $\A \x_1 x_k^\t$ and appears to be `more optimistic' than the real error, that is measured in the Frobenius norm.
In methods with the restricted pivoting  the Frobenius norm of $\A \x_1 x_{k+1}^\t \x_2 Y_l^\t \x_3 Z_m^\t$ is used to estimate the internal error and matches the real error more closely.

Timings for the approximation for different molecules and accuracy parameters are given in Tab.~\ref{tab2}. We compare all methods listed in Tab.~\ref{tab1} and the method proposed in~\cite{ost-chem-2009} based on incomplete cross approximation~\cite{tee-cross-2000} of unfoldings. For Wedderburn elimination methods we set $\pals=\ppow=3.$ 
From the preliminary experiments we see that this rather small number of `inner' iterations is sufficient; the experiments with $\pals$ and $\ppow$ set to $10,$ $30$ and $300$ show almost the same convergence and final accuracy for all the molecules.
For the reference we also provide time for one iteration of the Tucker-ALS method~\cite{tuckerals-1980,lathauwer-rank1-2000}, that uses output of Alg.~\ref{alg:wt3} as an initial guess.  
In our implementation of the Tucker-ALS the low-rank decomposition of~\eqref{eq:tals} is computed by the cross approximation, that is usually several ($2 \div 20$) times faster than the  SVD-based computations. 
However, even then one iteration of the  Tucker-ALS seems to be quite expensive. 
In practical computations several (usually $3 \div 20$)  iterations of the Tucker-ALS are required depending on the accuracy of the initial guess; it is also necessary to estimate the ranks of the desired approximation before starting the Tucker-ALS. 
In section~\ref{NUMs} we explain that each Tucker-ALS iteration requires $3r^2$ tenvecs, while Wedderburn methods require $\O(r)$ tenvecs to compute Tucker factors $U, V, W$ and $r^2$ tenvecs to generate the core $\G$ (see Tab.~\ref{tab1}). 
Therefore, for large ranks and fixed $\ppow$ and $\pals,$ Wedderburn methods become more efficient in comparison with a single iteration of the Tucker-ALS. 
We can note this behaviour in Tab.~\ref{tab2} on the lines corresponding to the high precision of approximation, where the ranks of the approximation are also high. 

We note that for this example the Wedderburn elimination with the restricted Lanczos-like pivoting (Alg.~\ref{alg:wt3}) is faster than all other versions of the Wedderburn elimination method. 
We also note that it outperforms the minimal Krylov recursion, since Alg.~\ref{alg:min} converges slowly (or even does not converge) and uses more iterations to reach the same accuracy level. 
WlncR also outperforms method based on cross approximation of unfoldings~\cite{ost-chem-2009}. 
We conclude that for this problem Alg.~\ref{alg:wt3} outperforms previously proposed methods.

\subsection{Recompression in operations with structured tensors}\label{NUMt2t}

The efficient operations with matrices and vectors in compressed tensor formats is crucial in the construction of efficient iterative methods for solving equations and eigenproblems in three and more dimensions. 
The approach to such highly-efficient tensor linear algebra subroutines was discussed in~\cite{ost-latensor-2009,st-chem-2009}, where it is shown that the efficient evaluation of all basic linear algebra subroutines with tensor-structured data is based on the fast recompression of certain structured tensor. As an example, consider Hadamard (elementwise) multiplication between $n_1 \x n_2 \x n_3$ tensors 
$$
\A = \G \x_1 U^{(A)} \x_2 V^{(A)} \x_3 W^{(A)}, \qquad
\B = \H \x_1 U^{(B)} \x_2 V^{(B)} \x_3 W^{(B)}.
$$
Let mode ranks of $\A$ be $r_1, r_2, r_3$ and mode ranks of $\B$ be $p_1,p_2,p_3.$ The result reads
\begin{equation}\label{eq:t2t}
\C = \F \x_1 U \x_2 V \x_3 W
\end{equation}
with $p_1r_1 \x p_2r_2 \x p_3r_3$ core $\F \eqdef \Kron(\G, \H)$ and \emph{non-orthogonal} factors  $U, V, W$ of sizes $n_1 \x p_1r_1,$ $n_2 \x p_2r_2 $ and $n_3 \x p_3r_3,$ respectively. Formally this is again the  Tucker format with the core and factors given by
$$
\F(ap, bq, cs) \eqdef \G(p,q,s) \H(a,b,c), \qquad U(i,ap) = U^{(A)}(i,p) U^{(B)}(i,a),
$$ 
and so on for $V, W.$
The mode ranks of $\C$ are products of correspondent mode ranks of $\A$ and $\B,$ and recompression is required to reduce the storage size.

In~\cite{st-chem-2009} the fast recompression method based on the individual filtering of factors was proposed. 
Numerical examples in~\cite{st-chem-2009} include evaluation of the Hartree potential for the electron density of molecules, discussed in section~\ref{NUMc}. 
This problem writes as a multiplication between three-level matrix given in the canonical format (with diagonal core tensor) and three-dimensional vector of the electron density given in the Tucker format. 
The computation requires $10 \div 120$ seconds depending on the complexity of molecule and desired accuracy of evaluation. However, examples of Tucker to Tucker multiplication were not presented in~\cite{st-chem-2009}, since this operation appears to be sufficiently more expensive. For the large molecules it requires up to an hour, that we consider not affordable.  

We show that fast and accurate multiplication between tensors given in the Tucker format can be done using Wedderburn-based methods.
Each tenvec operation with tensor~\eqref{eq:t2t} can be done in~$\O(q^4 + nq^2),$ where $q=\max(p,r),$ that is fast enough even for $n$ up to hundred thousands and  $p,r$ up to several hundreds.
We apply the discussed algorithm for the Hadamard multiplication of the discretized electron density of simple molecules to themselves. 
This operation can be a building block for algorithms that compute pointwise nonlinear functions of large tensors with linear in mode size complexity. 
One of the important applications is the cubic root of the electron density that appear in the Kohn-Sham model.
A good initial guess for such methods can be evaluated by the mimic algorithm~\cite{ost-chem-2009}. 

\begin{figure}[p]
\caption{Approximation accuracy of the Hadamard square of the methane electron density}
\begin{center}\hfil
\includegraphics[width=.45\textwidth]{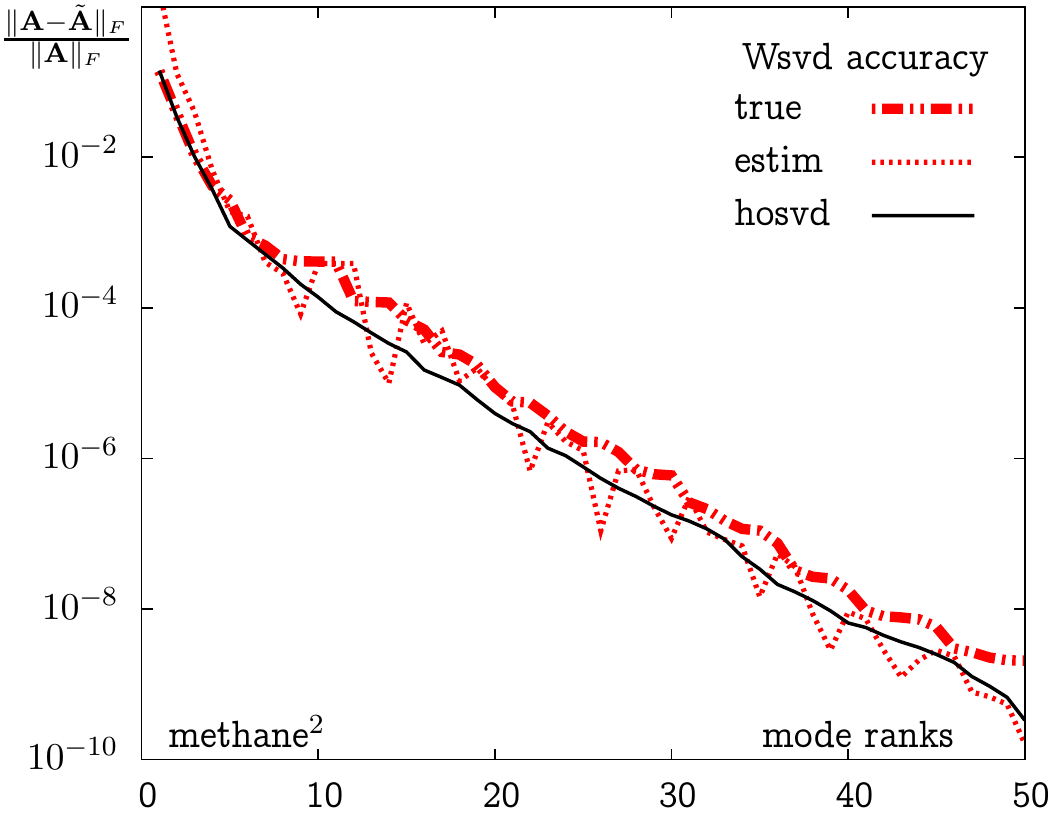} \hfil
\includegraphics[width=.45\textwidth]{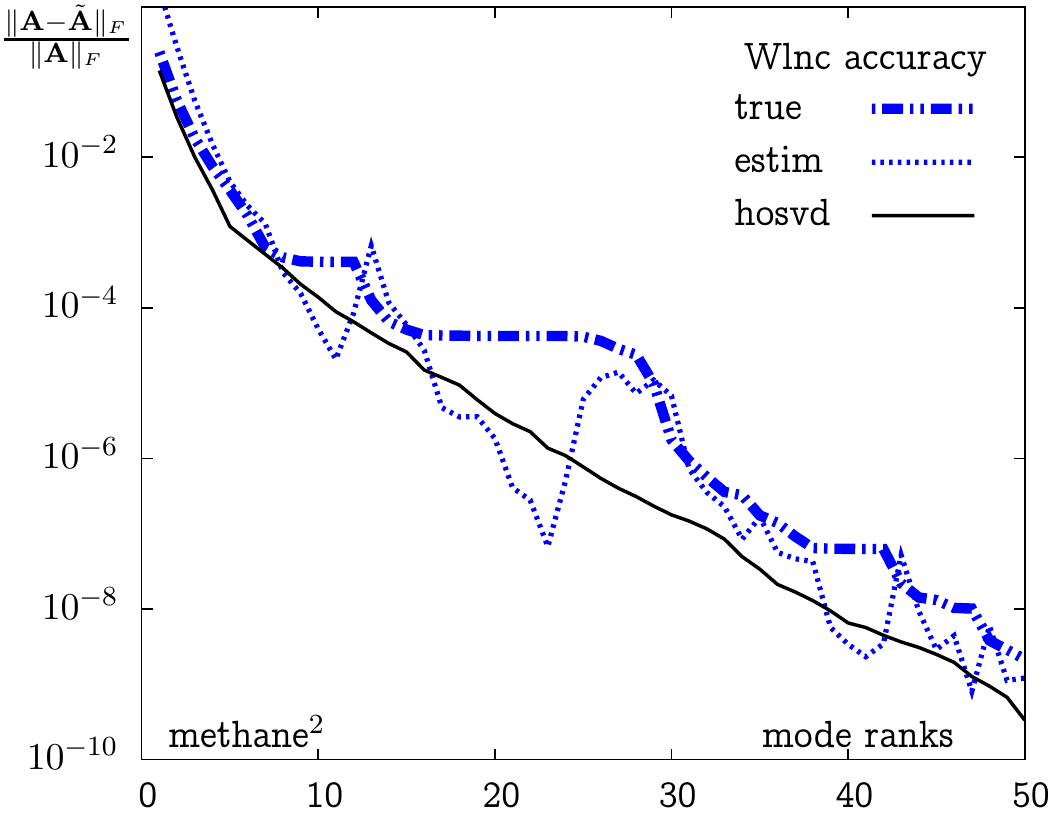} \hfil
\end{center}
\begin{center}\hfil
\includegraphics[width=.45\textwidth]{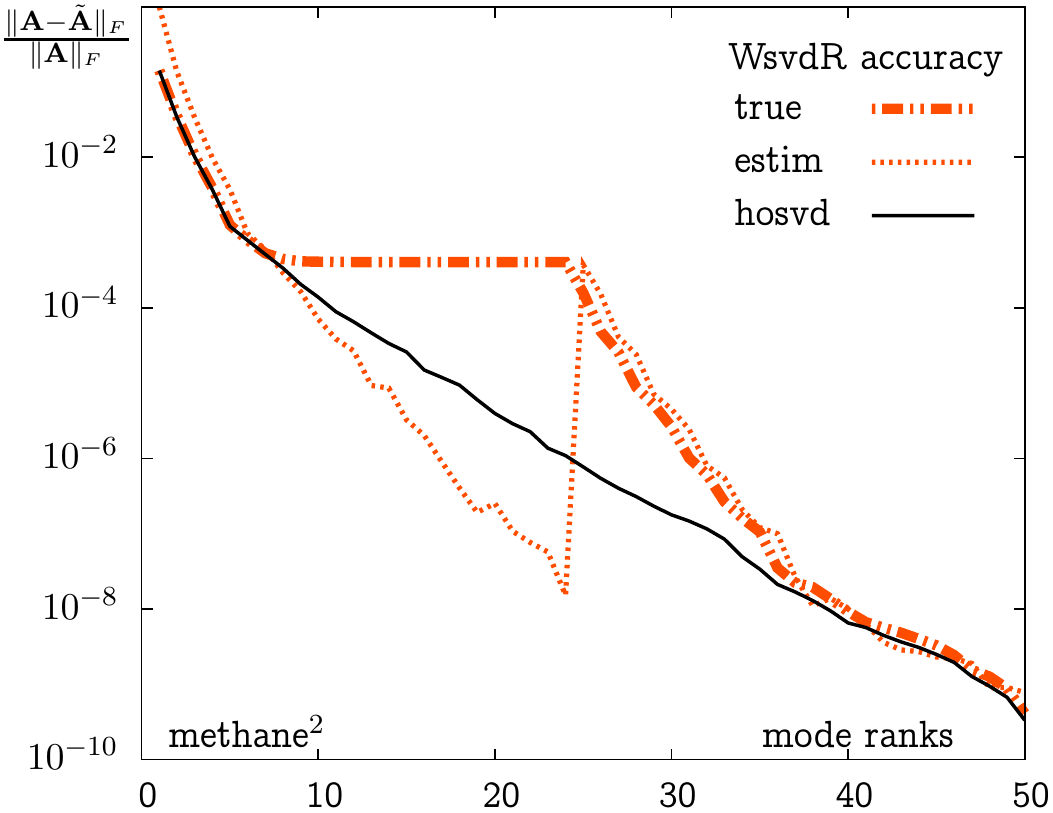} \hfil
\includegraphics[width=.45\textwidth]{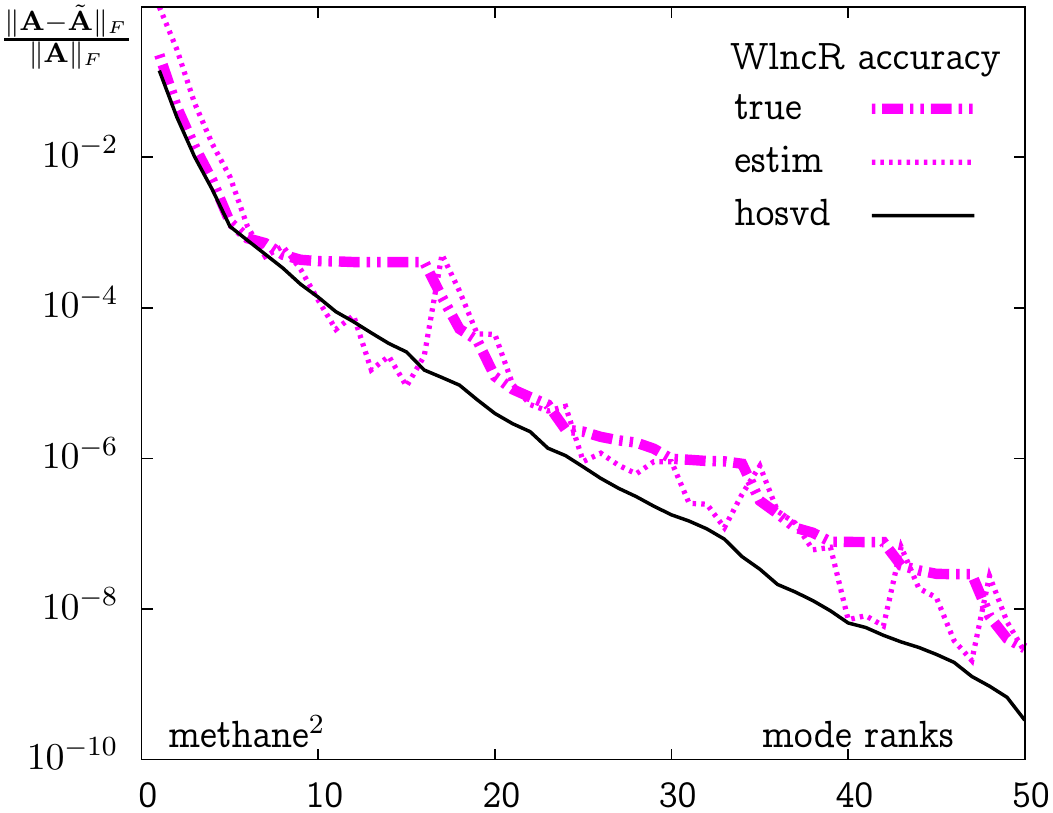} \hfil
\end{center}
\cntr{For the Hadamard square of the $5121 \x 5121 \x 5121$ tensor given in the Tucker format with mode ranks $(74,74,74),$ we compute Tucker approximation by algorithms, listed in Tab.~\ref{tab1}, and plot the relative accuracy in the Frobenius norm for different values of mode ranks. 
On each graph the thick dashed line shows the true error of approximation, the thin dashed line shows the estimate of error and the thin solid line shows the accuracy of the HOSVD approximation for the reference. 
} \label{fig:h1}
\end{figure}
\begin{figure}[p]
\caption{Approximation accuracy of the Hadamard square of the methane electron density}
\begin{center}\hfil
\includegraphics[width=.45\textwidth]{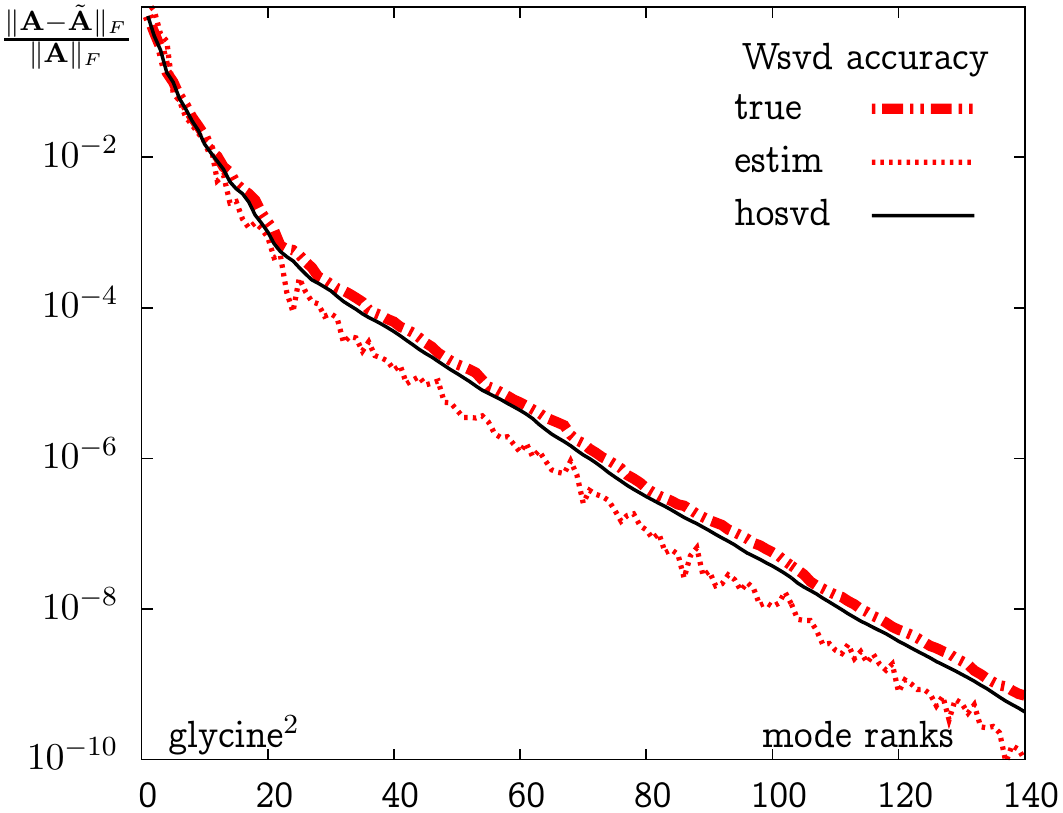} \hfil
\includegraphics[width=.45\textwidth]{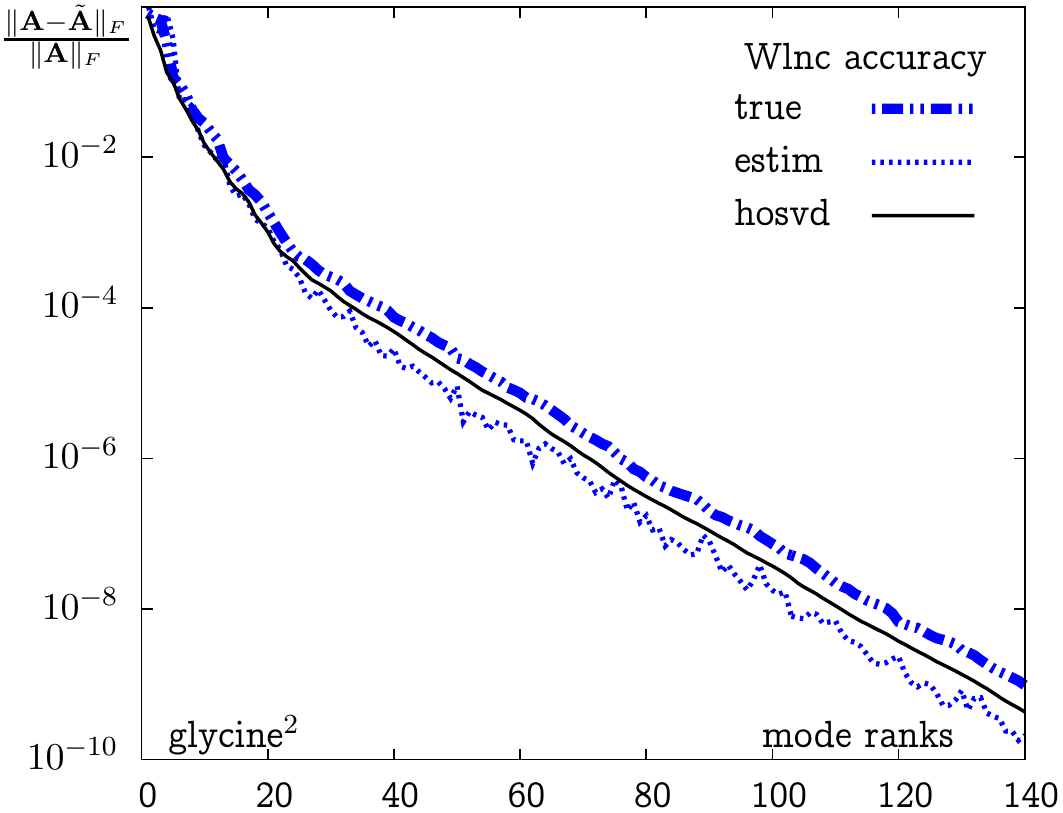} \hfil
\end{center}
\begin{center}\hfil
\includegraphics[width=.45\textwidth]{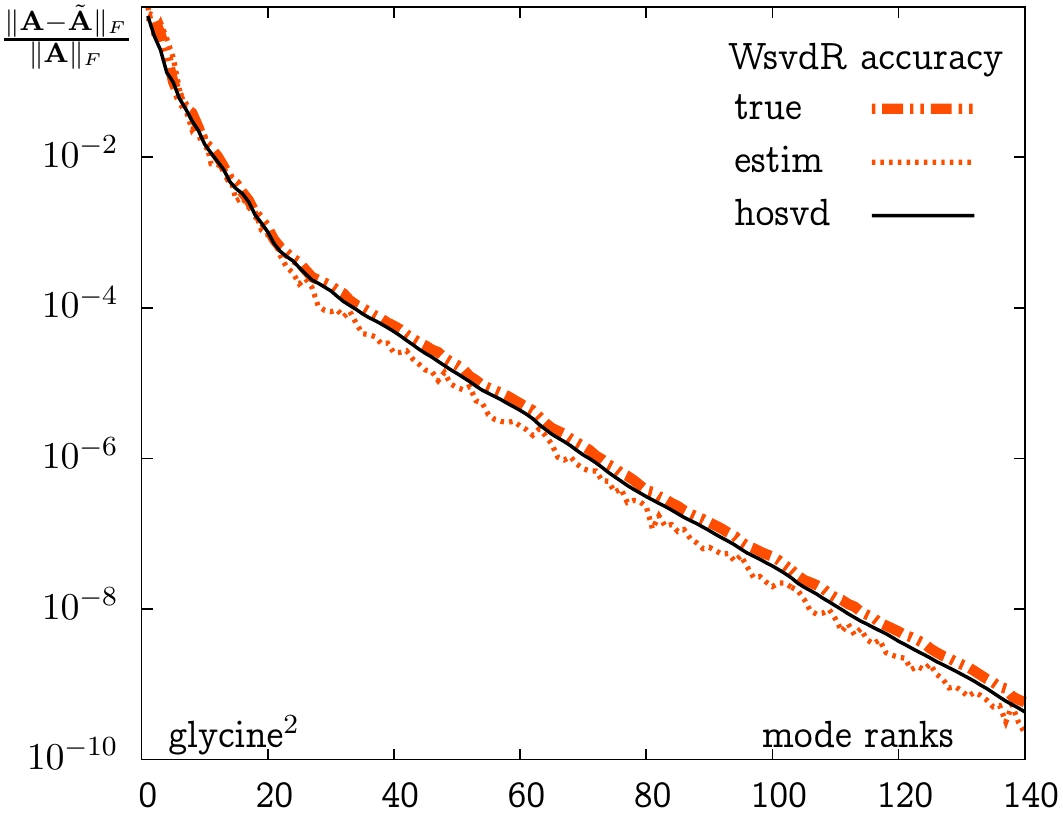} \hfil
\includegraphics[width=.45\textwidth]{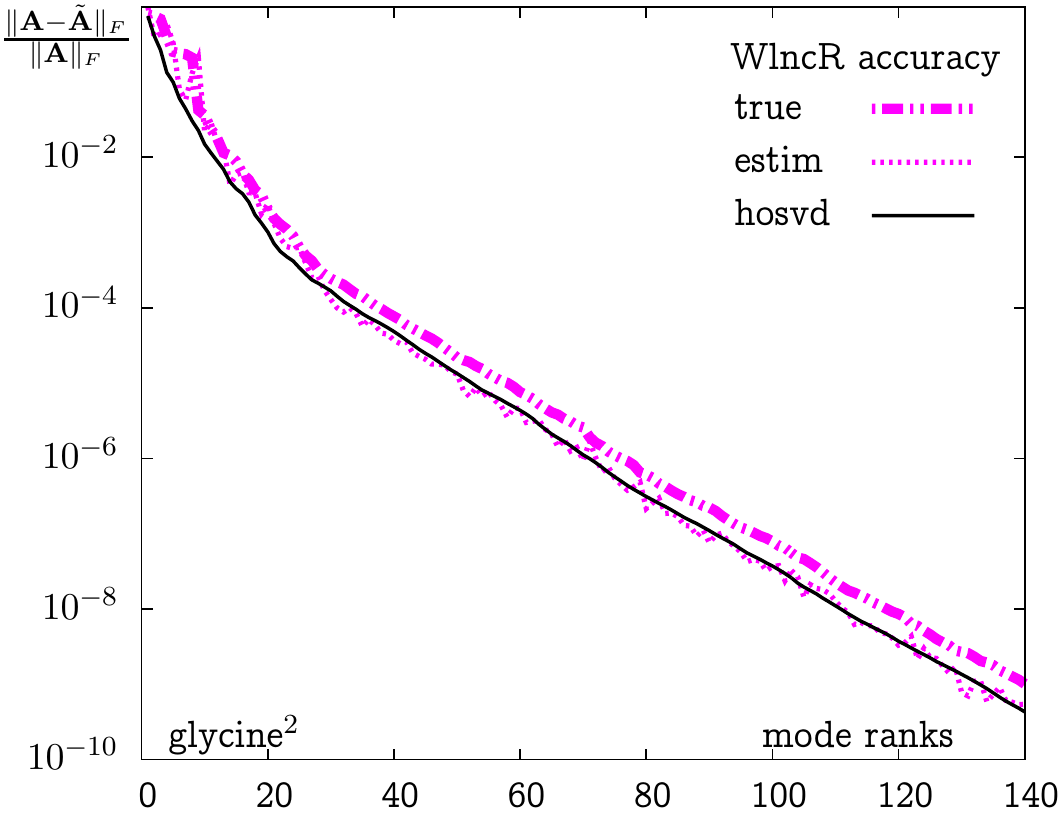} \hfil
\end{center}
\cntr{For the Hadamard square of the $5121 \x 5121 \x 5121$ tensor given in the Tucker format with mode ranks $(62,176,186),$ we compute Tucker approximation by algorithms, listed in Tab.~\ref{tab1}, and plot the relative accuracy in the Frobenius norm for different values of mode ranks. 
On each graph the thick dashed line shows the true error of approximation, the thin dashed line shows the estimate of error and the thin solid line shows the accuracy of the HOSVD approximation for the reference. 
} \label{fig:h2}
\end{figure}

The convergence of algorithms, i.e. the accuracy of rank-$(r,r,r)$ approximation for different $r$ is shown on Figs.~\ref{fig:h1},~\ref{fig:h2}. 
On each graph we compare accuracy of the  approximation computed by the Tucker-ALS~\cite{tuckerals-1980,lathauwer-rank1-2000} (thin solid line) with the accuracy of certain approximation method. 
The real accuracy $\|\A-\tilde\A\|_F$ was verified by comparing the result with the Tucker approximation computed by the Cross3D algorithm~\cite{ost-tucker-2008} with the accuracy set to $\eps=10^{-12}.$ The Cross3D algorithm was verified in~\cite{ost-tucker-2008,fkst-chem-2008} by  exhaustive check on cluster platforms.
The residual between two Tucker formats is computed as proposed in~\cite{ost-latensor-2009}.

As well as for recompression of electron density from canonical form, in this problem all versions of the Wedderburn elimination demonstrate good convergence, and for larger glycine molecule it is even more regular, than for methane. 
Comparing the different versions of Alg.~\ref{alg:wt1}, we note the same behaviour that is already described in section~\ref{NUMc}. 

Timings for the approximate computation of the Tucker-to-Tucker multiplication for different methods and accuracy parameters are given in Tab.~\ref{tab3}. 
We compare all methods listed in Tab.~\ref{tab1} and provide time for one iteration of the Tucker-ALS method~\cite{tuckerals-1980,lathauwer-rank1-2000}, that uses output of Alg.~\ref{alg:wt3} as an initial guess. For the Wedderburn elimination methods we set $\pals=\ppow=3.$
Again, we note that for this example the WlncR Alg.~\ref{alg:wt3} is faster than all other versions of Wedderburn elimination algorithms and also than one iteration of Tucker-ALS method. 
In this case the timings of different versions of Wedderburn elimination algorithms are more similar, because the cost of the evaluation of the core ($r^2$ tenvecs) dominates over $\O(nr^2)$ time for additional operations.
Nevertheless, we conclude that for this problem Alg.~\ref{alg:wt3} can be method of choice for fast approximate evaluation of operations with data in tensor formats.

\begin{table}[t]
\caption{Time for approximation of electron density}
\begin{center}
\begin{tabular}{cc|cccc|c}
molecule        & accuracy    & Wsvd   & Wlnc   & WsvdR  & WlncR  & TALS($1$)\\ \hline
methane         & $10^{-4}$   & $15.5$ & $13.8$ & $8.1$  & $2.8$  & $2.4$    \\
$(74,74,74)$    & $10^{-6}$   & $34$   & $33$   & $14.8$ & $9.7$  & $14.3$   \\
                & $10^{-8}$   & $63$   & $46$   & $43$   & $18.6$ & $42$     \\
                & $10^{-10}$  & $93$   & $74$   & $68$   & $40$   & $97$     \\ \hline
ethane          & $10^{-4}$   & $22$   & $20$   & $11.8$ & $4.2$  & $4.2$    \\
$(67,94,83)$    & $10^{-6}$   & $46$   & $41$   & $33$   & $14.0$ & $16.7$   \\
                & $10^{-8}$   & $82$   & $68$   & $72$   & $27$   & $45$     \\
                & $10^{-10}$  & $125$  & $105$  & $127$  & $56$   & $117$    \\ \hline
ethanol         & $10^{-4}$   & $120$  & $101$  & $106$  & $45$   & $52$     \\
$(128,127,134)$ & $10^{-6}$   & $281$  & $228$  & $293$  & $176$  & $257$    \\
                & $10^{-8}$   & $493$  & $419$  & $635$  & $441$  & $678$    \\
                & $10^{-10}$  & $736$  & $653$  & $1100$ & $808$  & $1370$   \\ \hline
glycine         & $10^{-4}$   & $179$  & $170$  & $177$  & $60$   & $64$    \\
$(62,176,186)$  & $10^{-6}$   & $442$  & $380$  & $600$  & $217$  & $270$   \\
                & $10^{-8}$   & $732$  & $600$  & $1033$ & $500$  & $646$   \\
                & $10^{-10}$  & $1010$ & $850$  & $1530$ & $888$  & $1223$  \\
\end{tabular}
\end{center}
\cntr{For the Hadamard square of the $5121 \x 5121 \x 5121$ tensor given in Tucker form, we compute Tucker approximation with different relative accuracy bound by different algorithms and show time in seconds. 
Algorithms Wsvd, Wlnc, WsvdR, WlncR are listed in Tab.~\ref{tab1}, 
`TALS(1)' is time for one iteration of Tucker-ALS method~\cite{tuckerals-1980,lathauwer-rank1-2000}
} \label{tab3}
\end{table}

\section{Conclusion and further work}
We presented the family of algorithms for the approximation of large three-dimensional tensor that access it only using the tenvec operation (tensor-by-vector-by-vector multiplication). 
Our approach is based on the Wedderburn rank-reduction formula and admits different strategies to select vectors of the Wedderburn elimination (`pivoting') that lead to algorithms with rather different convergence and  complexity estimates. 
The fastest algorithm from presented family, namely WlncR Alg.~\ref{alg:wt3}, may converge slowly or stagnate on the first steps of iterations. However, one can propose an efficient algorithm combining more efficient but demanding pivoting strategy on the first steps (for example, Wsvd, that is free from breakdowns) with the fast WlncR pivoting strategy on the next steps of algorithm.

The presented methods can be applied for the approximation of dominant subspaces of large structured tensors. 
In the provided numerical examples some of the proposed algorithms are much faster than the Tucker-ALS~\cite{tuckerals-1980,lathauwer-rank1-2000} algorithm and as fast as the minimal Krylov recursion~\cite{eldensavas-krylov-2009}, but more accurate in certain cases. 
The proposed methods can be directly applied for the fast computation of bilinear operations between structured tensors, but the efficiency can be further improved by combining them with the individual factor filtering proposed in~\cite{st-chem-2009}.

Canonical and Tucker formats can be straightforwardly generalized to $d$ dimensions, however each of them has serious drawbacks, and another decomposition should be used for high dimensions, for example recently introduced tensor-train (TT decomposition, see~\cite{ot-tt-2009,osel-compact-2009,osel-inside-2009}), that are based on the SVD techniques, but is free from the curse of dimensionality.
Therefore, it is very natural to extend the proposed ideas to the TT format.
In the sequel we will show how to construct algorithm for TT format that use tensor through tensor-by-vectors multiplication, i.e. Krylov-type methods in $d$ dimensions.

It is also important to analyse how does the choice of parameters $\ppow$ and $\pals$ change the convergence properties and final accuracy of the proposed methods. The preliminary experiments show that for the experiments provided in this paper it is sufficient to take $\pals=\ppow=3.$ The deeper theoretical analysis is required, that should generalize the theory of Arnoldi methods for the tensor case. This is a topic of a forthcoming work.

\section*{Acknowledgements}
Authors are grateful to Heinz-Juergen Flad and Rao Chinnamsettey for providing input data for the electron density function.
Authors are grateful to  Lars Eld\'en for inspiring paper~\cite{eldensavas-krylov-2009}, heartful attention to results of this paper, presented on ICSMT (Hong Kong, January 2010) and the $26^{\mathrm{th}}$ GAMM seminar (Leipzig, February 2010) and introducing us to the field of approximation of sparse tensors in network analysis.


\end{document}